\documentclass[twocolumn,secnumarabic,amssymb,
superscriptaddress,
nobibnotes,aps,
10pt]{revtex4-1}
\usepackage[toc,page,title,titletoc,header]{appendix}
\usepackage{stmaryrd}
\usepackage{amsmath,amsthm}
\usepackage{mathrsfs}
\usepackage{tabularx}
\usepackage{subfigure}
\usepackage{graphicx}
\usepackage{epsfig}
\usepackage{color}

\usepackage[colorlinks,linkcolor=blue, anchorcolor = blue, citecolor = blue, urlcolor = black]{hyperref}

\begin{document}

\title{A sharp-interface model for solid-state dewetting with wetting potential}

\author{Weijie Huang}
\affiliation{School of Mathematics and Statistics, Beijing Jiaotong University, Beijing 100044, China}
\email{wjhuang@bjtu.edu.cn}

\author{Xinran Ruan}
\affiliation{School of Mathematical Sciences, Capital Normal University, Beijing 100048, China}
\email{xinran.ruan@cun.edu.cn}


\begin{abstract}
We propose a sharp-interface model for solid-state dewetting of thin films with wetting potential, where the wetting effect is incorporated through a thickness-dependent surface energy. The model is governed by surface diffusion together with natural boundary conditions, and describes the morphological evolution of the film-vapor interface. For its numerical approximation, we develop an efficient semi-implicit finite element method based on a Taylor expansion of the wetting-potential term. Numerical simulations in two dimensions show that the proposed model and method can capture various dewetting phenomena. They also indicate that, as the range of the wetting potential tends to zero, the proposed model approaches the sharp-interface model with thickness-independent surface energy proposed in \cite{Wang15}. The model and numerical method are further extended to three dimensions, where the computations capture complex morphological evolution in solid-state dewetting.
\end{abstract}

\date{\today}
\maketitle


\section{Introduction}

Solid thin films deposited on substrates are often unstable in the as-deposited state. When heated to a temperature well below the melting point, a continuous film may agglomerate, break up into isolated parts, and eventually evolve into separated particles. This phenomenon is known as solid-state dewetting. It plays an important role in the morphological evolution of thin films and has attracted much attention in materials science, applied mathematics, and scientific computing.

\begin{figure}[htp]
\centering
\hspace*{-5mm}\includegraphics[width = .54\textwidth]{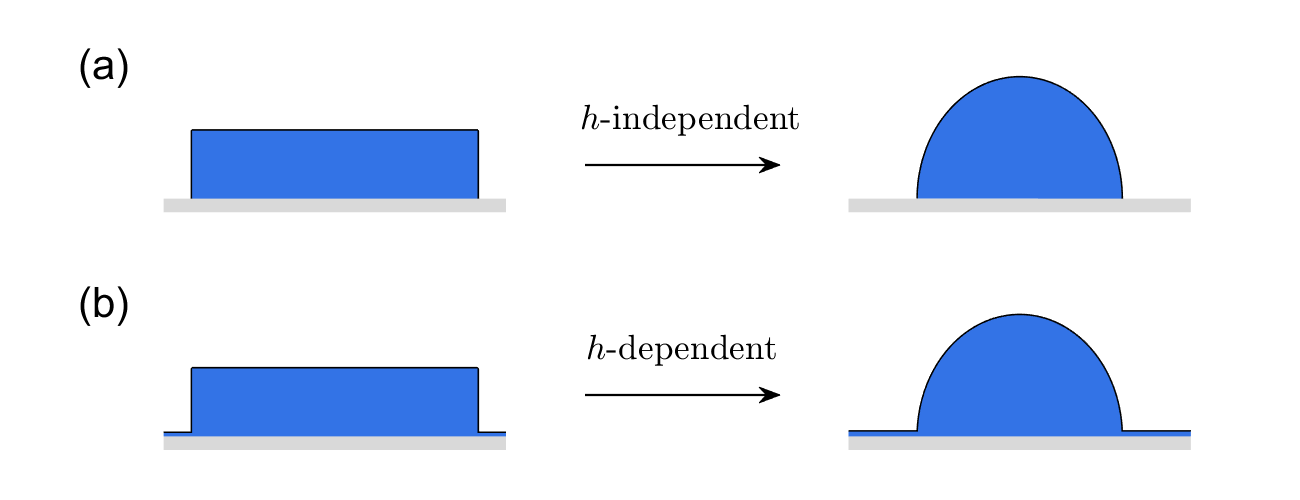}
\vspace*{-7mm}
\caption{A schematic illustration of the evolution from an initial film profile to its equilibrium shape for (a) \(h\)-independent surface energy and (b) \(h\)-dependent surface energy. In (b), a thin wetting layer remains on the substrate, and \(h\) denotes the local film thickness.}
\label{fig:dewetting0}
\end{figure}

From the modeling viewpoint, solid-state dewetting is driven by the minimization of the total surface free energy. In the classical setting, the film-vapor surface energy is assumed to be independent of the film thickness. In this case, dewetting leads to the shrinkage of the film and the exposure of the substrate, as illustrated in Fig.~\ref{fig:dewetting0}(a). The corresponding evolution is typically described by surface diffusion, supplemented with complicated boundary conditions at the moving contact line. Under this framework, several mathematical models and numerical methods have been developed for solid-state dewetting with thickness-independent surface energy, including sharp-interface models \cite{Srolovitz86a, Wong00, Wang15, Bao17a, Jiang20, Zhao20, BaoZhao23}, phase-field models \cite{Jiang12, Dziwnik17}, and kinetic Monte Carlo models \cite{Pierre09b}.

For very thin films, however, the film-vapor surface energy may be affected by the film-substrate interaction and thus depend on the local film thickness. The thickness-dependent part of the free energy is usually referred to as the wetting potential, and its derivative gives rise to the disjoining pressure. In this description, the substrate remains covered by a thin wetting layer even in the nominally dry region, and therefore no real moving contact line appears, as illustrated in Fig.~\ref{fig:dewetting0}(b). Wetting-potential effects in thin-film evolution have been studied in several different settings. Chiu et al. proposed a continuum boundary-layer model to account for the influence of the film-substrate interface on heteroepitaxial film growth \cite{Chiu95}. Peschka et al. studied liquid dewetting with intermolecular film-substrate interactions and reported rich dynamical behavior \cite{Peschka19}. For solid-state dewetting, Tripathi et al. derived triple-line kinetic boundary conditions in the presence of thickness-dependent wetting potential and mobility \cite{Tripathi18}, and later investigated the effect of disjoining pressure on mass shedding \cite{Tripathi20}. More recently, Zhou et al. proposed a regularized variational model for wetting/dewetting problems and studied its positivity-preserving property and asymptotic behavior, which further supports the use of wetting-potential-based descriptions \cite{Zhou24}. Nevertheless, for solid-state dewetting with wetting potential, sharp-interface dynamical models and efficient numerical approximations still deserve further investigation.

Motivated by these developments, we propose a sharp-interface model for solid-state dewetting with wetting potential and develop a semi-implicit finite element method for its numerical approximation. Owing to the presence of the wetting layer, the model avoids explicit moving contact lines and provides a convenient framework for describing complicated morphological evolution, including pinch-off, particle shedding, and the evolution of films with defects. We also investigate, in two dimensions, the relation between the present model and the sharp-interface model with thickness-independent surface energy \cite{Wang15, Bao17a} as the range of the wetting potential decreases. Furthermore, the model and numerical method are extended to three dimensions to simulate complex solid-state dewetting dynamics.

The rest of the paper is organized as follows. In Section 2, we present the sharp-interface model for solid-state dewetting with wetting potential in two dimensions. In Section 3, we introduce a semi-implicit finite element method for the proposed model. Section 4 is devoted to two-dimensional numerical results. In Section 5, we extend the model and the corresponding numerical method to the three-dimensional case. Numerical simulations in three dimensions are presented in Section 6. Finally, conclusions are drawn in Section 7.


\section{Wetting-potential model in two dimensions}

In this section, we present a sharp-interface model for solid-state dewetting with wetting potential in two dimensions.

\subsection{Surface free energy}

For a very thin solid film, the interaction between the film and the substrate cannot be neglected. As a result, the film-vapor surface energy depends on the film thickness when the film is only a few atomic layers thick. It is generally assumed that the film-vapor surface energy \(\gamma_{FV}:=\gamma(h)\) approaches the bulk value \(\gamma_\infty\) as the film thickness \(h\) increases. On the other hand, in the limit \(h\to 0\), the film-vapor surface energy is expected to approach \(\gamma_\infty\cos\theta_i\), where
\[
\theta_i=\arccos\Bigl(\frac{\gamma_{VS}-\gamma_{FS}}{\gamma_\infty}\Bigr)
\]
is the Young angle, and \(\gamma_{VS}\) and \(\gamma_{FS}\) are the vapor-substrate and film-substrate interfacial energies, respectively.

\begin{figure}[htp!]
\centering
\includegraphics[width=.45\textwidth]{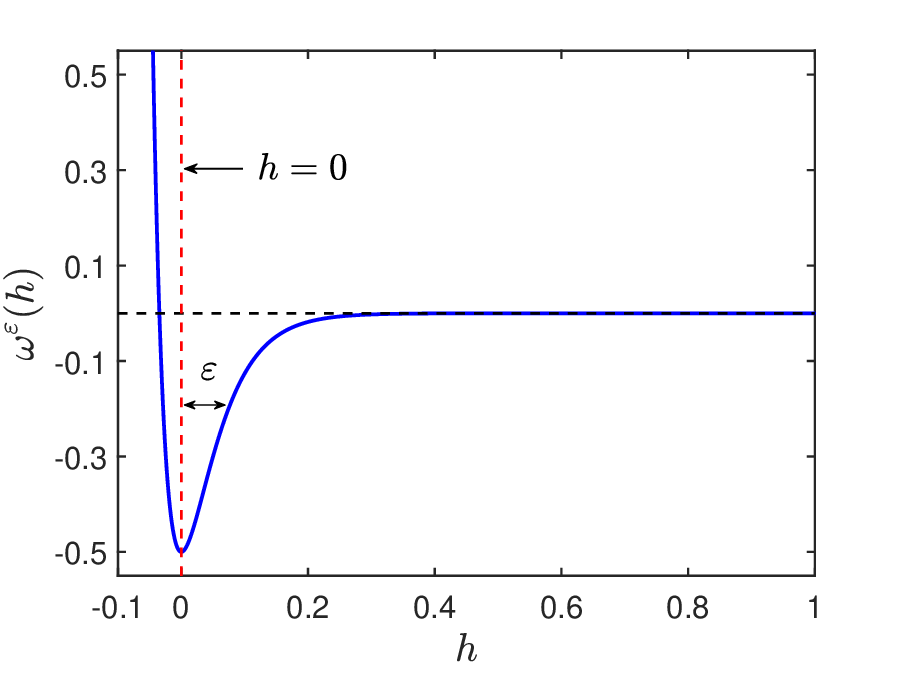}
\vspace*{-4mm}
\caption{Wetting potential as a function of the film thickness.}
\label{fig:wetting_potential}
\end{figure}

Based on these considerations, we choose the following exponentially decaying form for the dimensionless film-vapor surface energy
\begin{equation}\label{eq:surface_energy_gamma}
\gamma^\varepsilon(h) := 1 + (1 - \sigma)\bigl(e^{-\frac{h}{\varepsilon}} - 2e^{-\frac{h}{2\varepsilon}}\bigr),
\end{equation}
where \(\sigma=\cos\theta_i\), and \(\varepsilon\ll 1\) characterizes the range of the wetting interaction. We define
\[
\omega^\varepsilon(h):=\gamma^\varepsilon(h)-1
\]
as the wetting potential. Its profile is shown in Fig.~\ref{fig:wetting_potential}.

Moreover, for the choice \eqref{eq:surface_energy_gamma}, one has \((\gamma^\varepsilon)'(0)=0\), where the prime denotes differentiation with respect to \(h\), and \(\gamma^\varepsilon(h)\) attains its minimum at \(h=0\). Therefore, a thin wetting layer remains in the nominally dry region. Its thickness is determined by the parameter \(\varepsilon\), i.e., by the range of the wetting potential.

Finally, \(\gamma_\infty\omega^\varepsilon(0)=\gamma_{VS}-\gamma_{FS}-\gamma_\infty\), which coincides with the usual spreading coefficient. In this work, we focus on the partial wetting case, for which \(\omega^\varepsilon(0)<0\), as in \cite{Tripathi20}. In addition, since overhangs are excluded in the height-function formulation, all interface angles, including the equilibrium contact angle \(\theta_i\), are restricted to be smaller than \(\pi/2\).

\subsection{The model}

\begin{figure}[htp!]
\centering
\hspace{-4mm}
\includegraphics[width=.5\textwidth]{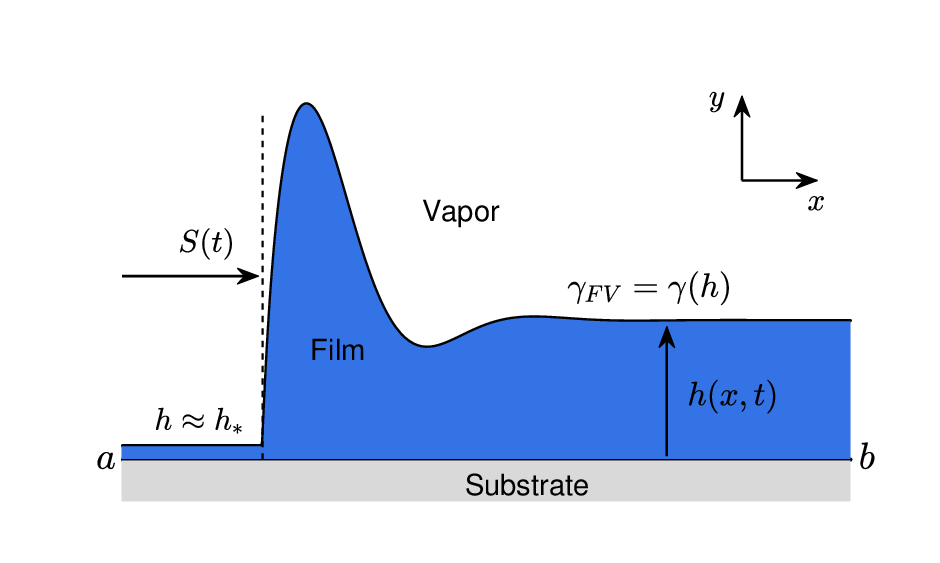}
\vspace*{-10mm}
\caption{A schematic illustration of the wetting-potential model for solid-state dewetting in 2D. Here \(h(x,t)\) denotes the height function, and a thin wetting layer with thickness \(h\approx h_*\) remains on the substrate in the nominally dry region.}
\label{fig:dewetting}
\end{figure}

We consider a thin solid film on a flat, rigid substrate in two dimensions, as illustrated in Fig.~\ref{fig:dewetting}. The total free energy of the system is given by
\begin{equation}\label{eq:energy_functional}
W^\varepsilon(h) = \int_\Gamma \gamma^\varepsilon(h)\, ds
= \int_a^b \gamma^\varepsilon(h) \sqrt{1 + (\partial_x h)^2}\, dx,
\end{equation}
where \(\Gamma=\Gamma(t)\) denotes the moving surface profile, namely the film-vapor interface, and \(s\) denotes the arc length along the interface.

The chemical potential is defined as the variational derivative of the free energy with respect to the height function \(h\), namely,
\begin{equation}\label{eq:chemical_potential}
\mu(x,t) := \dfrac{\delta W^{\varepsilon}(h)}{\delta h}
= \gamma^\varepsilon(h)\kappa + \dfrac{(\gamma^\varepsilon)'(h)}{\sqrt{1 + (\partial_x h)^2}}, \quad x\in I,
\end{equation}
where the first variation of the free energy is given in Appendix~\ref{Appendix:A}. Here
\[
\kappa = -\dfrac{\partial_{xx}h}{\bigl(1 + (\partial_x h)^2\bigr)^{3/2}}
\]
is the curvature of \(\Gamma\), and \(I=(a,b)\) with \(a<b\).

By surface diffusion kinetics \cite{Balluffi05}, the normal velocity \(V\) of the moving interface is given by
\begin{equation}
V = -\Omega_a \partial_s j, \qquad j = -m \partial_s \mu,
\end{equation}
where \(j\) is the surface mass flux, \(m\) is the surface mobility, assumed here to be constant, and \(\Omega_a\) is the volume per atom. Using the geometric relation
\begin{equation}
V = \dfrac{\partial_t h}{\sqrt{1 + (\partial_x h)^2}},
\end{equation}
together with
\[
\partial_s=\dfrac{1}{\sqrt{1+(\partial_x h)^2}}\partial_x,
\]
we obtain the evolution equation for the film height \(h\).

To obtain a dimensionless model, we further nondimensionalize the spatial variables and the parameter \(\varepsilon\) by a characteristic length \(R_0\). Choosing the time scale as \(R_0^4/(m\Omega\gamma_\infty)\), the dimensionless sharp-interface model for solid-state dewetting of a thin film on a flat substrate takes the form
\begin{equation}\label{eq:governing_equation}
\left\{
\begin{aligned}
& \partial_t h = \partial_x\Bigl(\dfrac{1}{\sqrt{1+(\partial_x h)^2}} \partial_x \mu\Bigr),\\
& \mu = \gamma^\varepsilon(h)\kappa + \dfrac{(\gamma^\varepsilon)'(h)}{\sqrt{1 + (\partial_x h)^2}},
\end{aligned}
\right.
\quad x\in I,\;\; t>0,
\end{equation}
where
\[
\kappa = -\dfrac{\partial_{xx}h}{\bigl(1 + (\partial_x h)^2\bigr)^{3/2}}.
\]
For simplicity, we continue to use the same symbols \(h\), \(t\), \(\mu\), \(\varepsilon\), and \(\kappa\) for the dimensionless variables.

The governing equation \eqref{eq:governing_equation} is supplemented with the zero-slope boundary conditions for \(h\) and the no-flux boundary conditions for \(\mu\):
\begin{gather}
\partial_x h(a,t) = 0, \qquad \partial_x h(b,t) = 0, \label{bc:h_x}\\
\partial_x \mu(a,t) = 0, \qquad \partial_x \mu(b,t) = 0. \label{bc:mass_zeros}
\end{gather}

The evolution equation \eqref{eq:governing_equation}, together with the boundary conditions \eqref{bc:h_x}--\eqref{bc:mass_zeros}, implies that the total film mass is conserved and the total free energy decreases monotonically during the evolution. See Appendix~\ref{Appendix:B} for details.

\section{Numerical method in two dimensions}

In this section, we derive a semi-implicit \(P_1\) finite element scheme for the sharp-interface model \eqref{eq:governing_equation}--\eqref{bc:mass_zeros}.

We first write the model in variational form. Given the initial data \(h(x,0)\), find \(h(\cdot,t)\in H^1(I)\) and \(\mu(\cdot,t)\in H^1(I)\) such that
\begin{equation}\label{eq:variational_formulation}
\left\{
\begin{aligned}
& (\partial_t h,\phi)
= -\Bigl(
\frac{\partial_x\mu}{\sqrt{1+(\partial_x h)^2}},\partial_x\phi
\Bigr),
\quad \forall\,\phi\in H^1(I),\\[4pt]
& (\mu,\psi)
=
\Bigl(
\frac{\gamma^\varepsilon(h)}{\sqrt{1+(\partial_x h)^2}}\partial_x h,\partial_x\psi
\Bigr) \\
& \qquad \quad +
\Bigl(
\dfrac{(\gamma^\varepsilon)'(h)}{\sqrt{1+(\partial_x h)^2}},\psi
\Bigr),
\quad \forall\,\psi\in H^1(I).
\end{aligned}
\right.
\end{equation}

Before introducing the fully discrete scheme, we first describe the treatment of the wetting term involving \((\gamma^\varepsilon)'(h)\). From \eqref{eq:surface_energy_gamma}, we have
\begin{equation}\label{eq:gamma_prime}
(\gamma^\varepsilon)'(h)
=
\frac{1-\sigma}{\varepsilon}
\bigl(-e^{-h/\varepsilon}+e^{-h/(2\varepsilon)}\bigr).
\end{equation}
For small \(h\), the nonlinear term \((\gamma^\varepsilon)'(h)\) is not suitable for a direct semi-implicit treatment. We therefore approximate \((\gamma^\varepsilon)'(h)\) near \(h=0\) by a quadratic function
\begin{equation}\label{eq:zeta_function}
\zeta(h)=c_1h+c_2h^2,
\end{equation}
where the coefficients \(c_1\) and \(c_2\) are chosen so that \(\zeta(h)\) and \((\gamma^\varepsilon)'(h)\) have the same value and first derivative at \(h=\bar h\ll1\). The construction of \(\zeta(h)\), together with the explicit expressions of \(c_1\) and \(c_2\), is given in Appendix~\ref{Appendix:C}. We then introduce the modified function
\begin{equation}\label{eq:gamma_prime_modified}
\widetilde{\gamma}'(h)=
\begin{cases}
\zeta(h), & h\le \bar h,\\[4pt]
(\gamma^\varepsilon)'(h), & h>\bar h.
\end{cases}
\end{equation}
This approximation allows us to construct a semi-implicit treatment of the wetting term in the numerical scheme.

Let \(I=[a,b]\) be partitioned uniformly into \(N\) subintervals \(I_j=[x_{j-1},x_j]\), \(j=1,2,\dots,N\), with mesh size
\[
\Delta x=\frac{b-a}{N},\qquad x_j=a+j\Delta x,\qquad j=0,1,\dots,N.
\]
We define the \(P_1\) finite element space by
\[
V_h:=\bigl\{v\in C(\bar I): v|_{I_j}\in P_1,\ j=1,2,\dots,N\bigr\}\subset H^1(I).
\]

Let \(\tau\) be the time step size and set \(t_m=m\tau\), \(m=0,1,2,\dots\). Denote by \(h^m\) and \(\mu^m\) the approximations of \(h(\cdot,t_m)\) and \(\mu(\cdot,t_m)\), respectively. Given \(h^m\in V_h\), we seek \(h^{m+1}\in V_h\) and \(\mu^{m+1}\in V_h\) such that
\begin{equation}\label{eq:scheme_1}
\Bigl( \frac{h^{m+1}-h^m}{\tau},\phi\Bigr)
+
\Bigl(
\frac{\partial_x\mu^{m+1}}{\sqrt{1+(\partial_x h^m)^2}},
\partial_x\phi
\Bigr)
=0,
\;\; \forall\,\phi\in V_h,
\end{equation}
and
\begin{equation}\label{eq:scheme_2}
\begin{aligned}
&(\mu^{m+1},\psi)
-
\Bigl(
\frac{\gamma^\varepsilon(h^m)\partial_x h^{m+1}}{\sqrt{1+(\partial_x h^m)^2}},
\partial_x\psi
\Bigr)  \\
&\qquad  =
\Bigl(
\dfrac{\widetilde{\gamma}'_{\mathrm{SI}}(h^m,h^{m+1})}{\sqrt{1+(\partial_x h^m)^2}},
\psi
\Bigr),
\;\;\; \forall\,\psi\in V_h,
\end{aligned}
\end{equation}
where \(\widetilde{\gamma}'_{\mathrm{SI}}(h^m,h^{m+1})\) is defined by
\begin{equation}\label{eq:gamma_prime_SI}
\widetilde{\gamma}'_{\mathrm{SI}}(h^m,h^{m+1})=
\begin{cases}
(c_1 +c_2 h^m)h^{m+1}, & h^m\le \bar h,\\[4pt]
(\gamma^\varepsilon)'(h^m), & h^m>\bar h.
\end{cases}
\end{equation}
Thus, when \(h^m\le \bar h\), the approximation is treated semi-implicitly, while for \(h^m>\bar h\), the original nonlinear term is evaluated explicitly at time level \(m\). As a result, only a linear system needs to be solved at each time step.

\section{Numerical simulations in two dimensions}
Unless otherwise stated, we choose the following stepped initial profile
\begin{equation}\label{eq:initial_data}
h(x,0)=\frac{1}{e^{-x+x_1}+1}+\frac{1}{e^{x-x_2}+1}-1,
\end{equation}
where \(x_1\) and \(x_2\) denote the positions of the two steps, with \(x_2>x_1\).

Since the migration of the triple line is driven by the wetting potential, a sufficiently fine mesh near the triple-line region is important in numerical simulations \cite{Tripathi18}. If the mesh is too coarse there, noticeable \(\varepsilon\)- and \(\sigma\)-dependent hysteresis effects may appear, which in turn affect the computed dynamics. Therefore, in all the following simulations, we use sufficiently fine grids with mesh size \(\Delta x\sim\varepsilon\). Unless otherwise stated, we also take \(\bar h=\varepsilon\).

\subsection{Equilibrium shape of a small island}

\begin{figure}[htp!]
\centering
\hspace*{-5mm}\includegraphics[width = .52\textwidth]{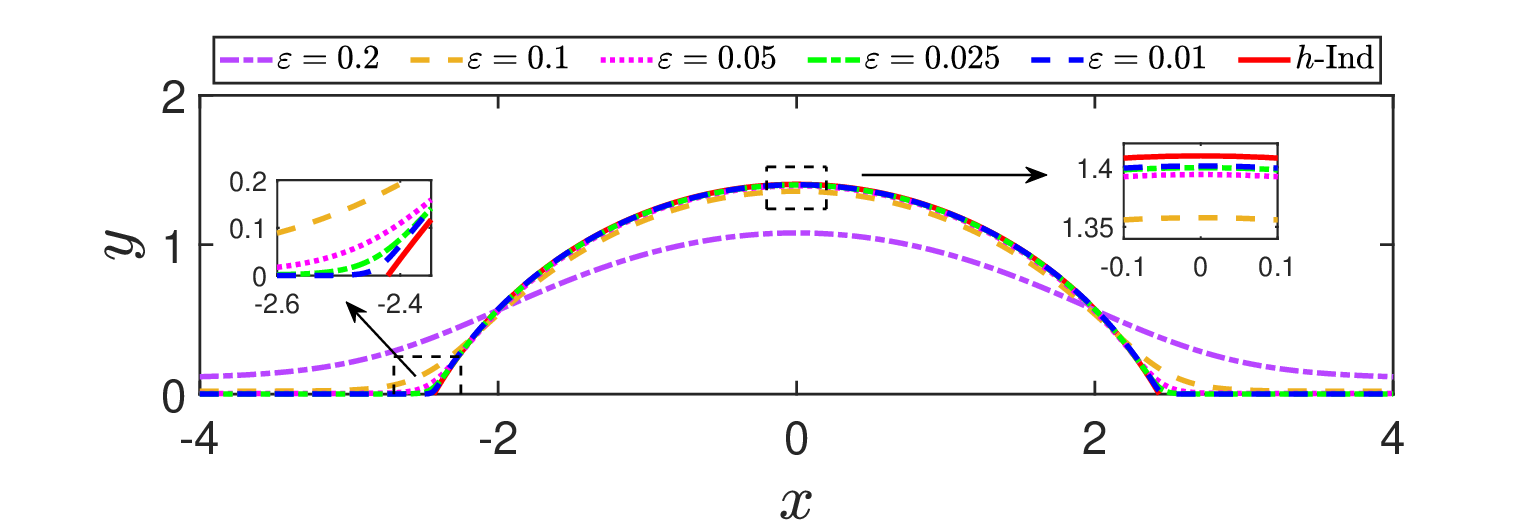}
\vspace*{-7mm}
\caption{Comparison of the numerical equilibrium shapes (dashed lines) of a small island with the theoretical equilibrium shape (solid red line) of the \(h\)-independent model \cite{Bao17b}, where \(\sigma=\cos(\pi/3)\).}
\label{fig:convergence_small_island}
\end{figure}

In this subsection, we study the equilibrium shape of a small island for several values of \(\varepsilon\) while fixing \(\sigma\). Starting from the initial profile \eqref{eq:initial_data} with \(x_1=-2.5\) and \(x_2=2.5\), Fig.~\ref{fig:convergence_small_island} shows that the computed equilibrium shapes approach the theoretical equilibrium shape of the \(h\)-independent model predicted by the generalized Winterbottom construction \cite{Bao17b} as \(\varepsilon\) decreases.

\begin{figure}[htp!]
\centering
\includegraphics[width = .42\textwidth]{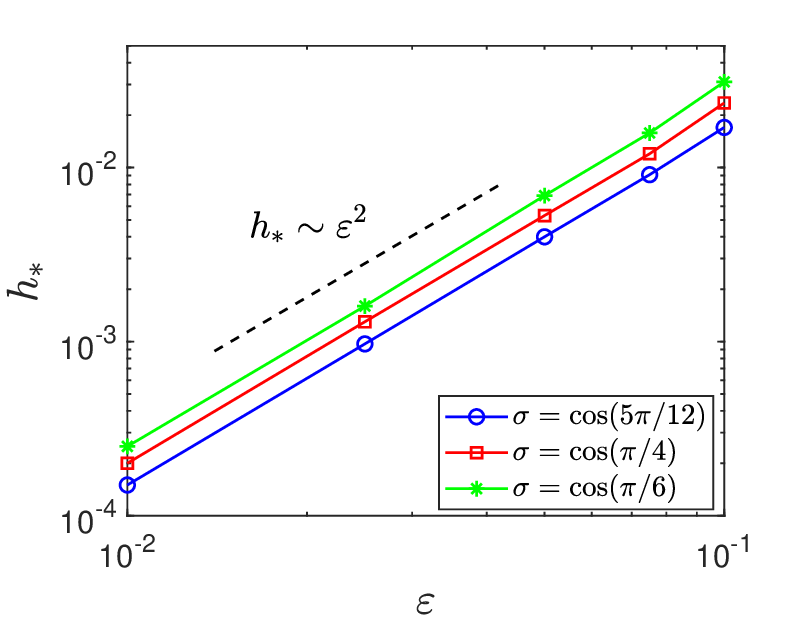}
\vspace*{-3mm}
\caption{Thickness of the wetting layer in equilibrium.}
\label{fig:thickness_wetting_layer}
\end{figure}

We also examine the thickness of the wetting layer in equilibrium. We denote by \(h_*\) the nearly uniform film thickness away from the island. As shown in Fig.~\ref{fig:thickness_wetting_layer}, \(h_*\) decreases quadratically as \(\varepsilon\) decreases, which is consistent with the asymptotic result \(h_*=O(\varepsilon^2)\) for the precursor thickness in the regularized model \cite{Zhou24}. Moreover, \(h_*\) decreases as \(\sigma\) decreases, that is, as the Young's angle \(\theta_i\) increases.

\subsection{Dynamics of finite islands}

\begin{figure}[htp!]
\centering
\includegraphics[width = .45\textwidth]{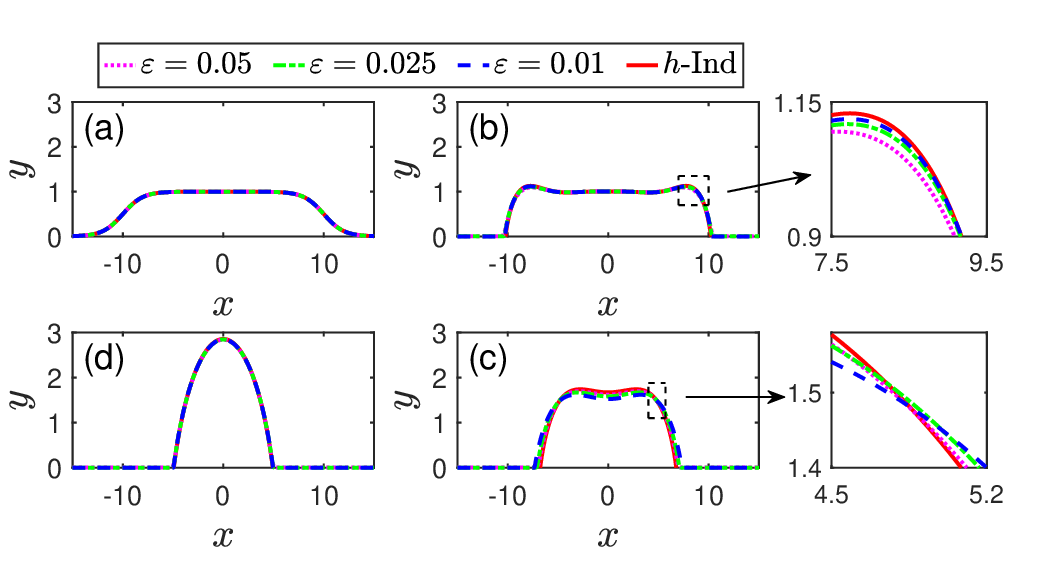}
\vspace*{-3mm}
\caption{Snapshots of the evolution of an initially small island under the proposed model for several values of \(\varepsilon\), together with the corresponding result of the \(h\)-independent model \cite{Bao17a}: (a) \(t=0\), (b) \(t=5\), (c) \(t=50\), (d) \(t=1000\). Here \(\sigma=\cos(\pi/3)\) and \(\eta=100\), where \(\eta\) is the contact line mobility in the \(h\)-independent model. Note the different vertical and horizontal scales.}
\label{fig:evolution_small_island}
\end{figure}

We first consider the evolution of a small island. Starting from the initial profile \eqref{eq:initial_data} with \(x_1=-10\) and \(x_2=10\), Fig.~\ref{fig:evolution_small_island} shows several snapshots of the evolution obtained by the proposed method for several values of \(\varepsilon\), with \(\sigma=\cos(\pi/3)\). The corresponding profiles are close to those obtained from the \(h\)-independent model solved by the parametric finite element method (PFEM) \cite{Bao17a}, and the agreement improves as \(\varepsilon\) decreases.

\begin{figure}[htp!]
\centering
\includegraphics[width = .45\textwidth]{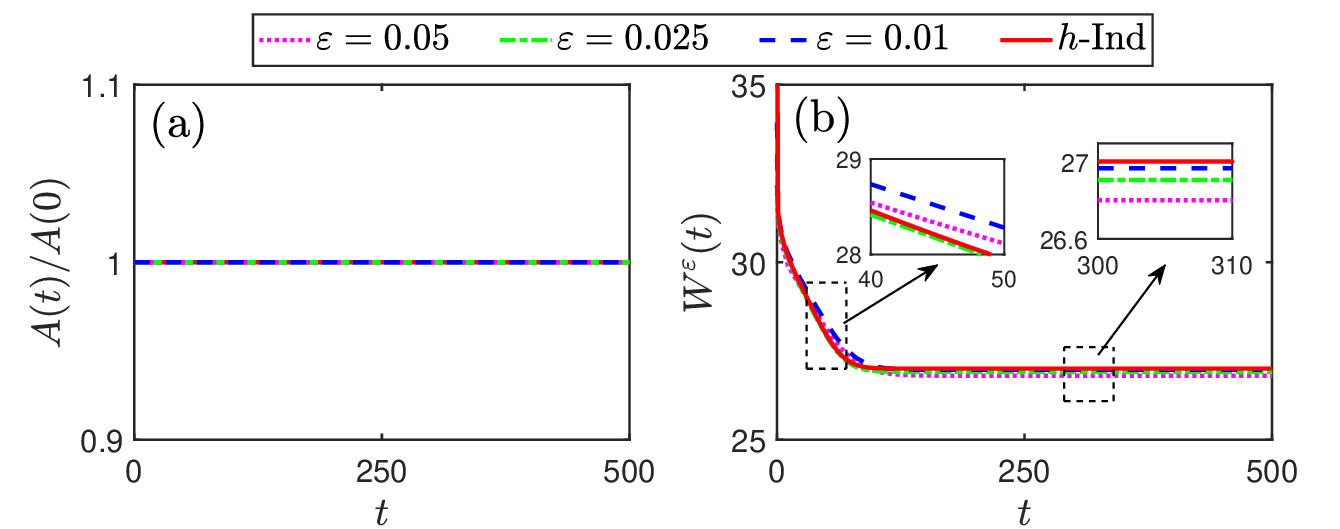}
\vspace*{-3mm}
\caption{Evolution of (a) the normalized area \(A(t)/A(0)\) and (b) the total free energy \(W^\varepsilon(t)\) defined by \eqref{eq:energy_functional} for the example in Fig.~\ref{fig:evolution_small_island}.}
\label{fig:small_numerical_property}
\end{figure}

Figure~\ref{fig:small_numerical_property} shows the evolution of the normalized area \(A(t)/A(0)\) and the total free energy \(W^\varepsilon(t)\). The total mass is well conserved throughout the evolution, while the free energy decreases monotonically in time. Moreover, as \(\varepsilon\) decreases, the free-energy curves become closer, especially near equilibrium.

\begin{figure}[htp!]
\centering
\hspace*{-0mm}\includegraphics[width = .45\textwidth]{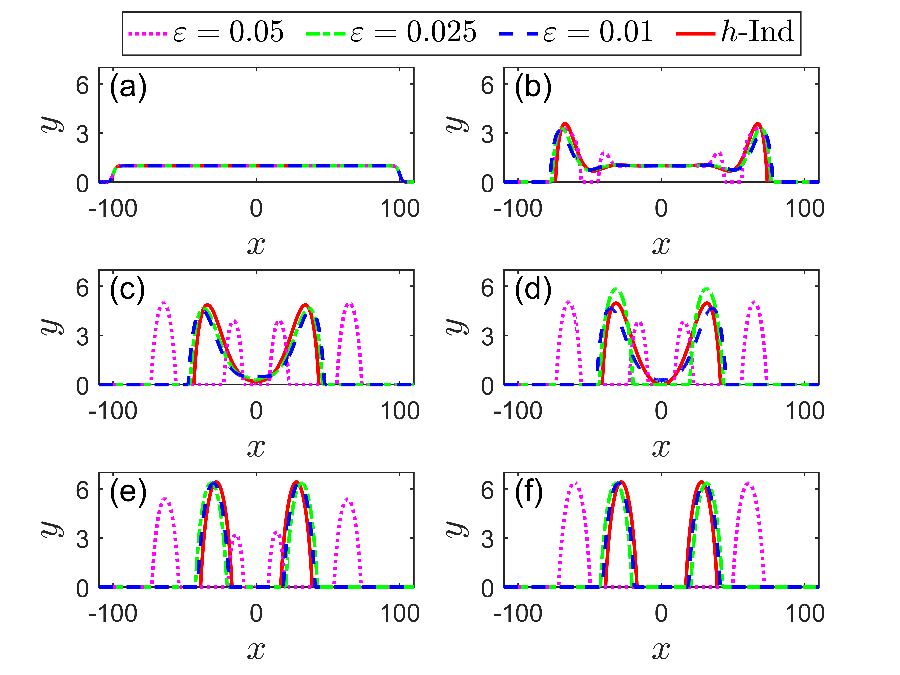}
\vspace*{-3mm}
\caption{Snapshots of the evolution of an initially long island under the proposed model for several values of \(\varepsilon\), together with the corresponding result of the \(h\)-independent model: (a) \(t=0\), (b) \(t=1000\), (c) \(t=2500\), (d) \(t=14710\), (e) \(t=16340\), (f) \(t=20000\). Here \(\sigma=\cos(\pi/3)\) and \(\eta=100\). Note the different vertical and horizontal scales.}
\label{fig:evolution_large_island}
\end{figure}

We next consider a long island by taking \(x_1=-100\) and \(x_2=100\) in \eqref{eq:initial_data}. It is well known that pinch-off may occur when the aspect ratio of the initial island is sufficiently large \cite{Dornel06,Jiang12,Wang15}. Figure~\ref{fig:evolution_large_island} shows the evolution for several values of \(\varepsilon\), with \(\sigma=\cos(\pi/3)\), where pinch-off is clearly observed. Again, the numerical results obtained by the proposed model approach those of the \(h\)-independent model as \(\varepsilon\) decreases.

A notable advantage of the present model is that pinch-off is captured naturally, without any explicit treatment of topological changes. By contrast, in the \(h\)-independent model, such changes must be handled manually. In our PFEM computations, whenever the valley depth drops below \(10^{-5}\), the lowest point is moved to the substrate and the film is split into two parts at the new contact point.
We also observe that a larger \(\varepsilon\) leads to faster valley thinning and earlier contact with the substrate. As a result, when \(\varepsilon=0.05\), the long island pinches off into four particles, whereas for \(\varepsilon=0.025\) and \(0.01\), only two particles remain. This acceleration of mass shedding induced by the disjoining pressure is consistent with the observations in \cite{Tripathi20}. Moreover, for \(\varepsilon=0.05\), the smaller particles gradually shrink and are eventually absorbed by the two larger ones, as shown in Fig.~\ref{fig:evolution_large_island}(d)--\ref{fig:evolution_large_island}(f). This coarsening behavior was also reported in \cite{Tripathi20}.

\begin{figure}[htp!]
\centering
\includegraphics[width = .45\textwidth]{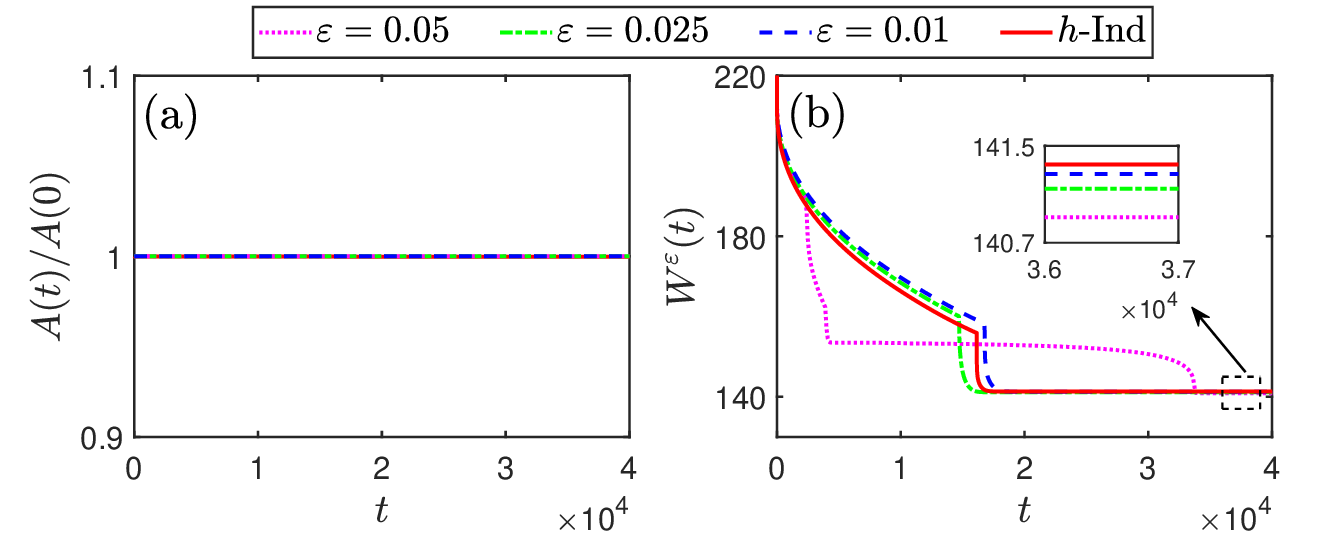}
\vspace*{-3mm}
\caption{Evolution of (a) the normalized area \(A(t)/A(0)\) and (b) the total free energy \(W^\varepsilon(t)\) defined by \eqref{eq:energy_functional} for the example in Fig.~\ref{fig:evolution_large_island}.}
\label{fig:large_numerical_property}
\end{figure}

Figure~\ref{fig:large_numerical_property} shows the corresponding evolution of the normalized area and the total free energy. Mass is again well conserved throughout the evolution. The free energy decreases monotonically in time, with sharp drops occurring immediately after pinch-off and when the small particles are completely absorbed.

\begin{figure}[htp!]
\centering
\hspace*{-2mm}\includegraphics[width = .49\textwidth]{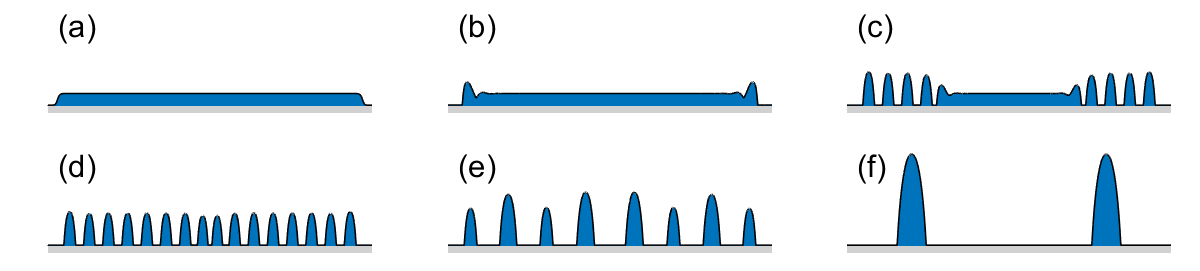}
\vspace*{-7mm}
\caption{Snapshots of the evolution of a long island film with aspect ratio \(300\) under the proposed model: (a) \(t=0\), (b) \(t=100\), (c) \(t=500\), (d) \(t=1000\), (e) \(t=1\times10^{4}\), (f) \(t=2\times10^{5}\). Here \(\sigma=\cos(\pi/3)\) and \(\varepsilon=0.1\). Note the different vertical and horizontal scales.}
\label{fig:large_eat_small_01}
\end{figure}

\begin{figure}[htp!]
\centering
\hspace*{-2mm}\includegraphics[width = .49\textwidth]{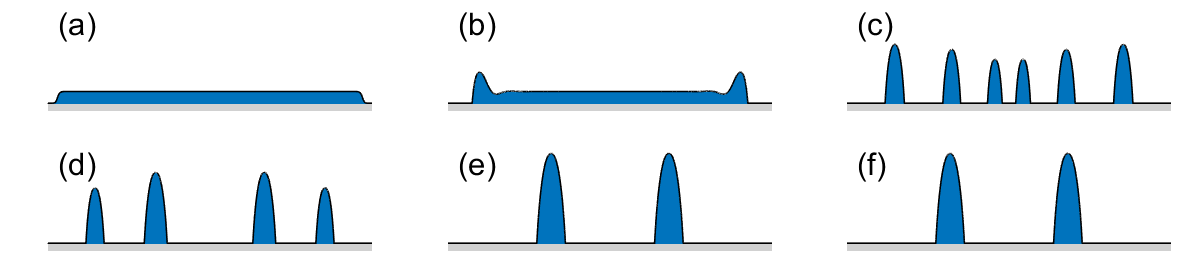}
\vspace*{-7mm}
\caption{Snapshots of the evolution of a long island film with aspect ratio \(300\) under the proposed model: (a) \(t=0\), (b) \(t=1\times10^{3}\), (c) \(t=1\times10^{4}\), (d) \(t=5\times10^{4}\), (e) \(t=1\times10^{5}\), (f) \(t=2\times10^{5}\). Here \(\sigma=\cos(\pi/3)\) and \(\varepsilon=0.025\). Note the different vertical and horizontal scales.}
\label{fig:large_eat_small_0025}
\end{figure}

\begin{figure}[htp!]
\centering
\hspace*{-2mm}\includegraphics[width = .49\textwidth]{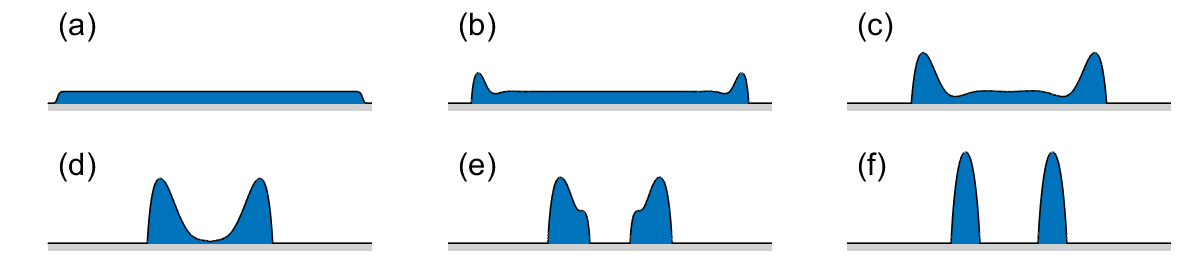}
\vspace*{-7mm}
\caption{Snapshots of the evolution of a long island film with aspect ratio \(300\) under the proposed model: (a) \(t=0\), (b) \(t=1\times10^{3}\), (c) \(t=2\times10^{4}\), (d) \(t=5.739\times10^{4}\), (e) \(t=5.77\times10^{4}\), (f) \(t=1\times10^{5}\). Here \(\sigma=\cos(\pi/3)\) and \(\varepsilon=0.01\). Note the different vertical and horizontal scales.}
\label{fig:large_eat_small_001}
\end{figure}

\begin{figure}[htp!]
\centering
\hspace*{-2mm}\includegraphics[width = .49\textwidth]{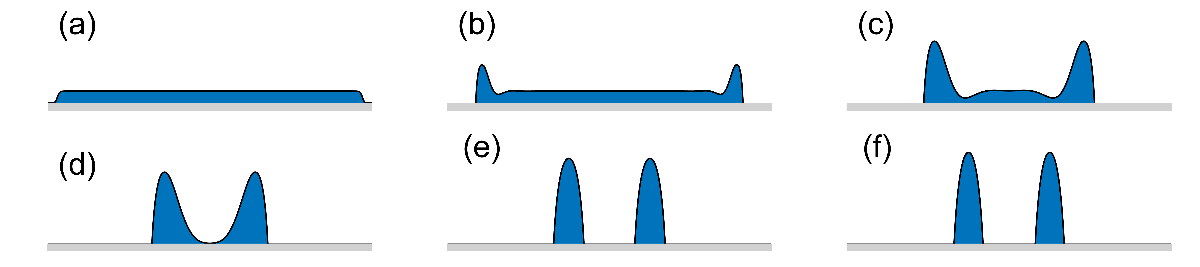}
\vspace*{-7mm}
\caption{Snapshots of the evolution of a long island film with aspect ratio \(300\) obtained by the \(h\)-independent model: (a) \(t=0\), (b) \(t=1\times10^{3}\), (c) \(t=2\times10^{4}\), (d) \(t=4.588\times10^{4}\), (e) \(t=4.8\times10^{4}\), (f) \(t=1\times10^{5}\). Here \(\sigma=\cos(\pi/3)\) and \(\eta=100\). Note the different vertical and horizontal scales.}
\label{fig:large_eat_small_PFEM}
\end{figure}

To further examine the influence of \(\varepsilon\) on coarsening, we consider a longer island with \(x_1=-150\) and \(x_2=150\) in \eqref{eq:initial_data}. Figures~\ref{fig:large_eat_small_01}--\ref{fig:large_eat_small_PFEM} compare the morphological evolution predicted by the proposed model for several values of \(\varepsilon\) with that of the \(h\)-independent model. As \(\varepsilon\) decreases, the coarsening effect becomes much weaker, and the evolution for \(\varepsilon=0.01\) is already very close to that of the \(h\)-independent model. The remaining difference in the pinch-off time between Fig.~\ref{fig:large_eat_small_001} and Fig.~\ref{fig:large_eat_small_PFEM} is mainly due to the choice of the contact line mobility \(\eta\).

\subsection{Number of agglomerates formed from finite islands}

\begin{figure}[htp!]
\centering
\includegraphics[width = .43\textwidth]{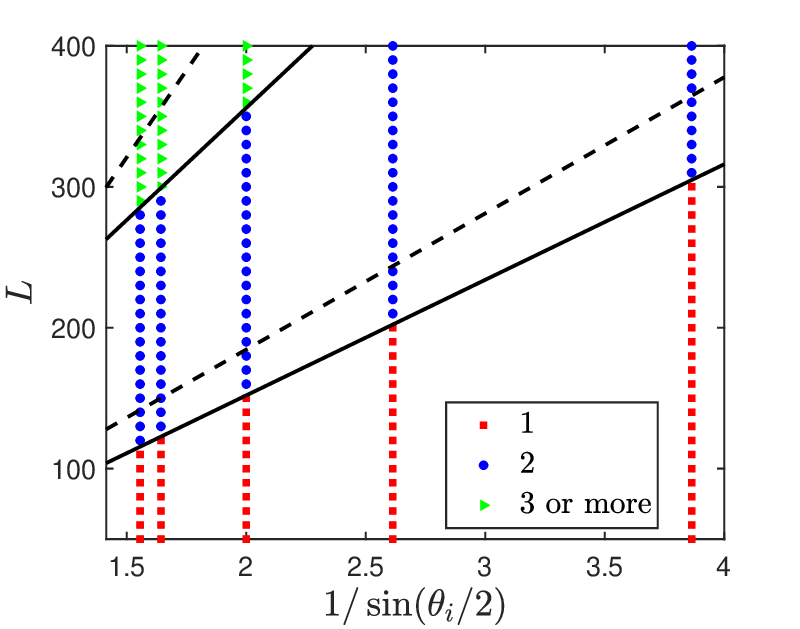}
\vspace*{-2mm}
\caption{Number of agglomerates formed from a high-aspect-ratio island versus the initial length \(L=x_2-x_1\) and the Young angle \(\theta_i\) for \(\varepsilon=0.01\). The dashed black line shows the numerical result of Dornel \cite{Dornel06} for rectangular islands with length \(L\) and height \(1\).}
\label{fig:pinch_off_numbers}
\end{figure}

As discussed above, when the aspect ratio is sufficiently large, a thin film may undergo pinch-off and eventually break into several agglomerates. The number of agglomerates depends on both the Young angle \(\theta_i\) and the parameter \(\varepsilon\), as illustrated in Fig.~\ref{fig:evolution_large_island}. To examine this dependence more systematically, we fix \(\varepsilon=0.01\) and perform a series of numerical simulations for large islands with different aspect ratios and different values of \(\theta_i\). The results are summarized in Fig.~\ref{fig:pinch_off_numbers}. We observe clear boundaries separating the regions with 1, 2, and 3 (or more) agglomerates. For comparison, the corresponding results for the \(h\)-independent model reported by Dornel \cite{Dornel06} are also shown.

\subsection{Dynamics of semi-infinite films}

\begin{figure}[htp!]
\centering
\hspace*{-5mm}\includegraphics[width = .54\textwidth]{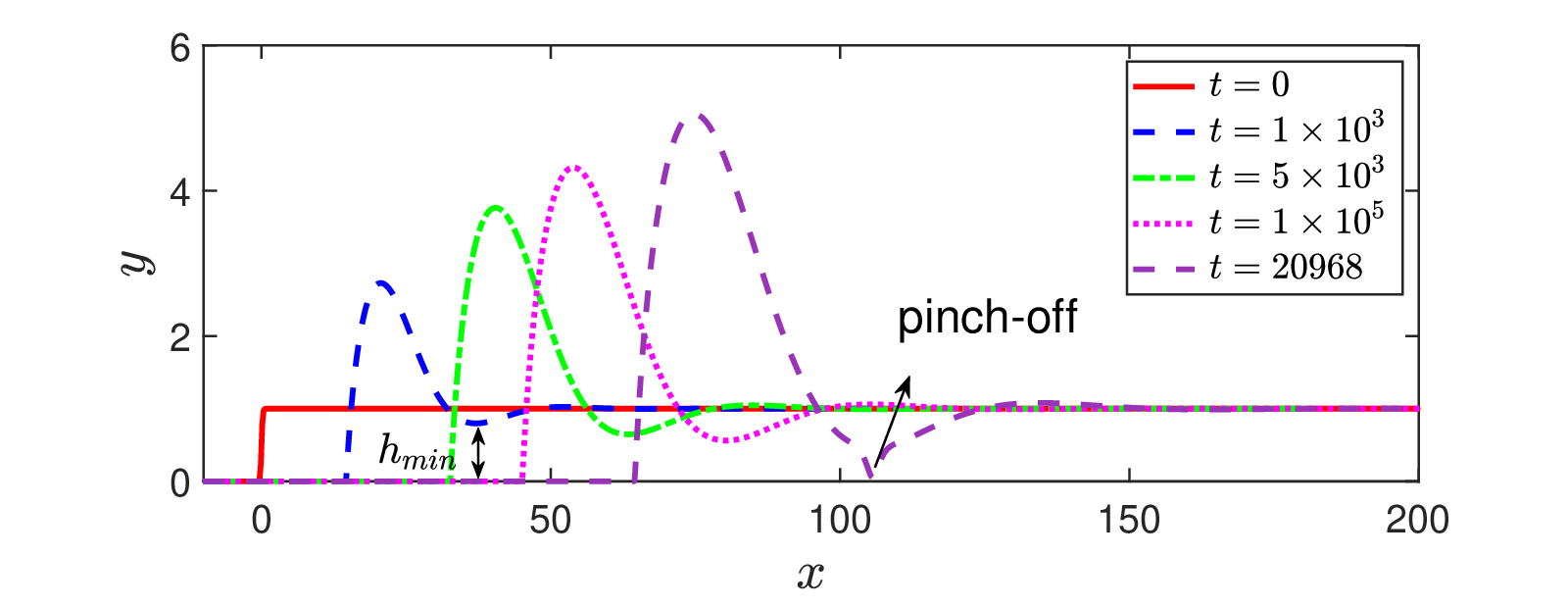}
\vspace*{-7mm}
\caption{Evolution of a semi-infinite film with \(\varepsilon=0.025\) and \(\sigma=\cos(\pi/3)\). Note the different vertical and horizontal scales.}\label{fig:evolution_semi_infty}
\end{figure}

In this subsection, we consider the evolution of initially semi-infinite films described by \eqref{eq:initial_data} with \(x_1=0\) and \(x_2\) chosen sufficiently large, for example \(x_2=10^5\). Figure~\ref{fig:evolution_semi_infty} shows a typical example with \(\sigma=\cos(\pi/3)\) and \(\varepsilon=0.025\). At the far-field end, we still impose the natural boundary condition and the zero-flux condition, while allowing the film height to vary. Whenever \(|h(b)-1|\ge 10^{-6}\), the computational domain is extended to the right by one unit so that the new right endpoint becomes \(b+1\), and the corresponding mesh points are added. As shown in Fig.~\ref{fig:evolution_semi_infty}, solid-state dewetting first leads to the formation of ridges near the film edge, followed by the development of a valley. As time evolves, the valley initially sinks slowly. Once its height falls below a certain threshold, it approaches the substrate very rapidly due to the wetting potential.

\begin{figure}[htp!]
\centering
\includegraphics[width = .43\textwidth]{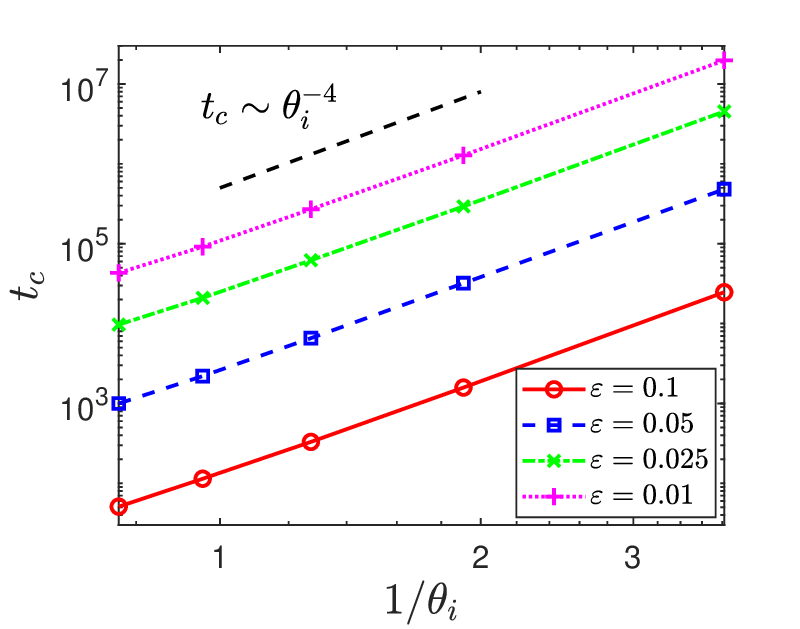}
\vspace*{-2mm}
\caption{First mass-shedding time \(t_c\) as a function of the Young angle \(\theta_i\).}\label{fig:time_first_pinchoff}
\end{figure}

To quantify this process, we first record the time \(t_c\) of the first mass-shedding event for several values of \(\varepsilon\) and \(\theta_i\). Figure~\ref{fig:time_first_pinchoff} plots \(t_c\) as a function of the Young angle \(\theta_i\), and indicates the scaling law
$$t_c\sim\theta_i^{-4},$$
which agrees well with the theoretical result for the \(h\)-independent model under the small-slope assumption obtained by Wong \cite{Wong00}. In addition, \(t_c\) increases significantly as \(\varepsilon\) decreases, but it does not show a simple power-law dependence on \(\varepsilon\), which is also consistent with the observation in \cite{Tripathi20}.

\begin{figure}[htp!]
\centering
\includegraphics[width = .43\textwidth]{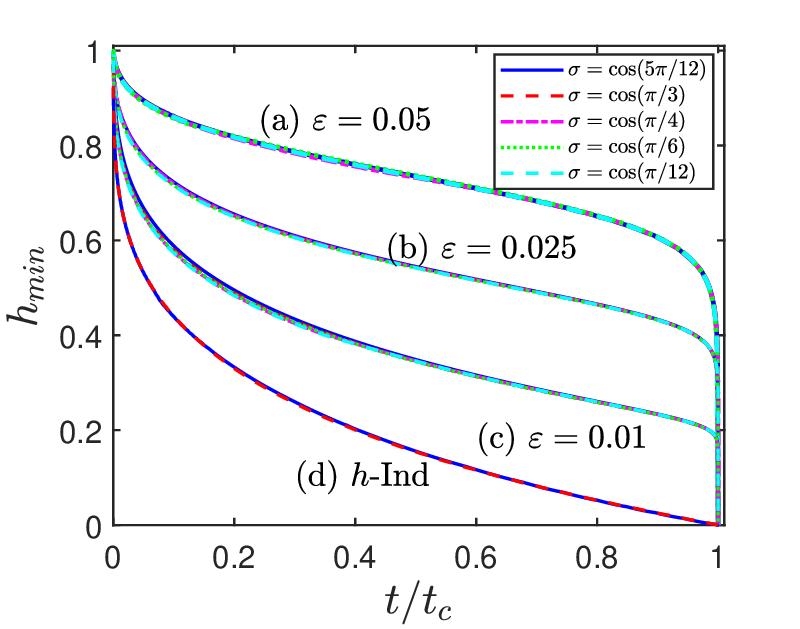}
\vspace{-3mm}
\caption{Evolution of the minimum valley thickness \(h_{\min}(t)\) before the first mass-shedding event. The small discrepancies between curves with the same \(\varepsilon\) and different \(\sigma\) are attributed to insufficient mesh resolution and to the fact that \(\bar h\) is not taken sufficiently small.}
\label{fig:h_min}
\end{figure}

We next define \(h_{\min}(t)\) as the minimum valley thickness of the film profile in the wetting region. Figure~\ref{fig:h_min} shows the evolution of \(h_{\min}\) as a function of \(t/t_c\) before the first mass-shedding event. Ignoring some minor differences, we observe that the evolution of \(h_{\min}\) is essentially independent of \(\sigma\). Moreover, \(h_{\min}\) drops sharply once it falls below a certain critical value, and this critical value is positively correlated with \(\varepsilon\), that is, with the range of the wetting potential.

\begin{figure}[htp!]
\centering
\hspace*{-2mm}\includegraphics[width = .5\textwidth]{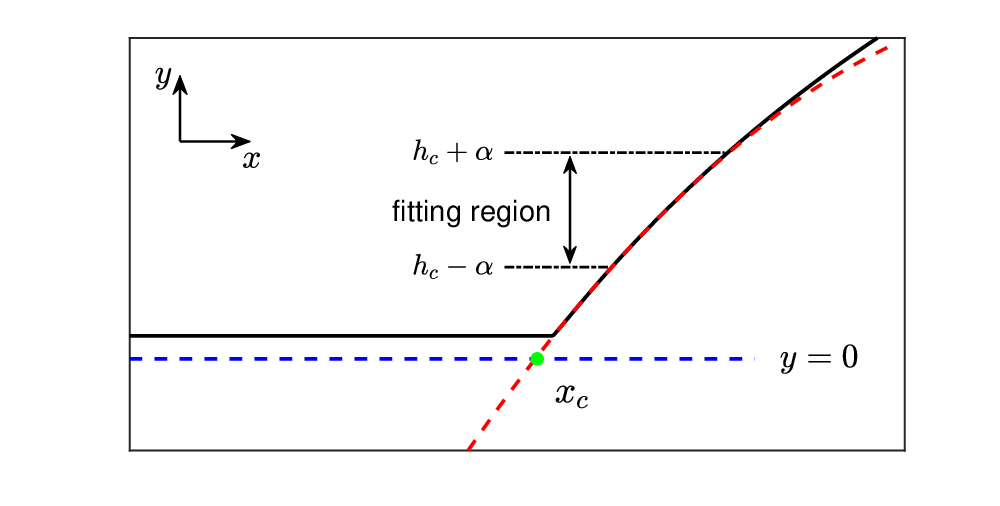}
\vspace*{-7mm}
\caption{Schematic of the fitting procedure used to define the effective contact point position \(x_c\). The solid black curve denotes the film profile, the dashed red curve denotes the quadratic fitting in the prescribed region, and the dashed blue line denotes the substrate.} \label{fig:contact_point_position}
\end{figure}

Finally, we examine the motion of an effective contact point. In the present \(h\)-dependent model, there is no genuine contact point in the usual sharp-interface sense. To facilitate comparison with the \(h\)-independent model, we therefore introduce a contact point through the fitting procedure shown in Fig.~\ref{fig:contact_point_position}. More precisely, it is defined as the intersection of the substrate \(y=0\) and a quadratic fitting of the film profile in a prescribed region. In the following computations, we choose \(h_c=0.2\) and \(\alpha=0.1\).

\begin{figure}[htp!]
\centering
\hspace*{-5mm}\includegraphics[width = .51\textwidth]{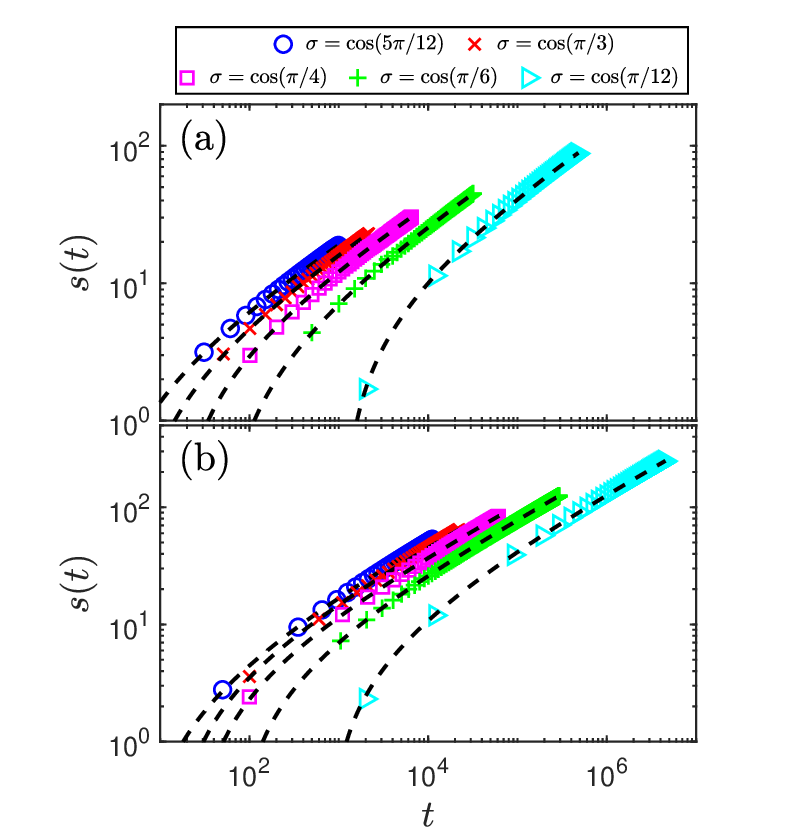}
\vspace*{-5mm}
\caption{Effective contact point position versus time during the first mass-shedding cycle, displayed on log--log axes, for (a) \(\varepsilon=0.05\) and (b) \(\varepsilon=0.025\). The dashed black lines show the fitting function \eqref{eq:fitting_contact_point}.}
\label{fig:fitting_contact_line}
\end{figure}

We next fit the position of the effective contact point by
\begin{equation}\label{eq:fitting_contact_point}
S(t)=c+at^{0.4}+bt^{0.2}.
\end{equation}
This fitting law was originally derived for the Mullins model with \(h\)-independent isotropic film--vapor surface energy under a fixed contact angle and vanishing mass-flux boundary conditions \cite{Wong00}. Although the present \(h\)-dependent model does not contain a genuine contact point in the sharp-interface sense, Fig.~\ref{fig:fitting_contact_line} shows that the numerically defined effective contact point still agrees well with \eqref{eq:fitting_contact_point} over a wide time interval. In particular, for both \(\varepsilon=0.05\) and \(\varepsilon=0.025\), the curves corresponding to different values of \(\sigma\) follow the same overall trend predicted by this law. This suggests that, during the first mass-shedding cycle, the motion of the effective contact point in the present model is still well captured by the fitting law from the \(h\)-independent theory.

\section{Wetting potential model in 3D}
As in the two-dimensional case, the sharp-interface model for three-dimensional solid-state dewetting can be written as (see Appendix~\ref{Appendix:D} for details)
\begin{equation}\label{eq:governing_equation_3D}
\left\{
\begin{aligned}
& h_t = \sqrt{1 + |\nabla h|^2}\,\Delta_s \mu,\\
& \mu = \gamma^\varepsilon(h)\kappa + \dfrac{(\gamma^\varepsilon)'(h)}{\sqrt{1 + |\nabla h|^2}},
\end{aligned}
\right.
\qquad (x,y)\in\Omega,\quad t>0,
\end{equation}
where \(\Delta_s=\nabla_s\cdot\nabla_s\) is the surface Laplace operator, \(\nabla_s\) is the surface gradient, \(\Omega=(a,b)\times(c,d)\), and \(\kappa\) denotes the curvature of the film-vapor interface \(S\). The corresponding boundary conditions are
\begin{gather}
\nabla h\cdot \mathbf{n} = 0, \quad \text{on } \partial\Omega, \\
\nabla_s \mu \cdot \mathbf{c} = 0, \quad \text{on } \partial S,
\end{gather}
where \(\mathbf{n}\in\mathbb{R}^2\) is the outward unit normal vector on \(\partial\Omega\), and \(\mathbf{c}\in\mathbb{R}^3\) is the outward unit conormal vector along \(\partial S\).

The corresponding variational formulation and semi-implicit finite element discretization follow the same strategy as in the two-dimensional case, and are therefore omitted here. In the next section, we present several three-dimensional numerical examples to illustrate the capability of the wetting potential model in capturing typical dewetting morphologies and their dependence on the parameter \(\varepsilon\).

\section{Numerical simulations in three dimensions}

In this section, we present several three-dimensional numerical examples to demonstrate the capability of the wetting potential model in capturing typical dewetting morphologies. We consider three classes of initial geometries: square islands, elongated cuboid islands, and more complex-shaped islands. Most of the corresponding examples have also been studied in the \(h\)-independent setting \cite{Jiang20,Zhao20}. In those studies, however, the computations are generally not continued after topological changes take place.

Unless otherwise specified, we take
\[
\sigma=\cos(4\pi/9), \qquad \varepsilon=0.05, \qquad \bar h=\varepsilon.
\]

\subsection{Evolution of square islands with increasing size}
We first investigate square islands of different sizes to illustrate the transition from edge retraction and corner accumulation to hole formation and eventual breakup. The initial wetting-layer thickness is set to \(10^{-5}\).

\begin{figure}[htp!]
\centering
\hspace*{-3mm}\includegraphics[width=.49\textwidth]{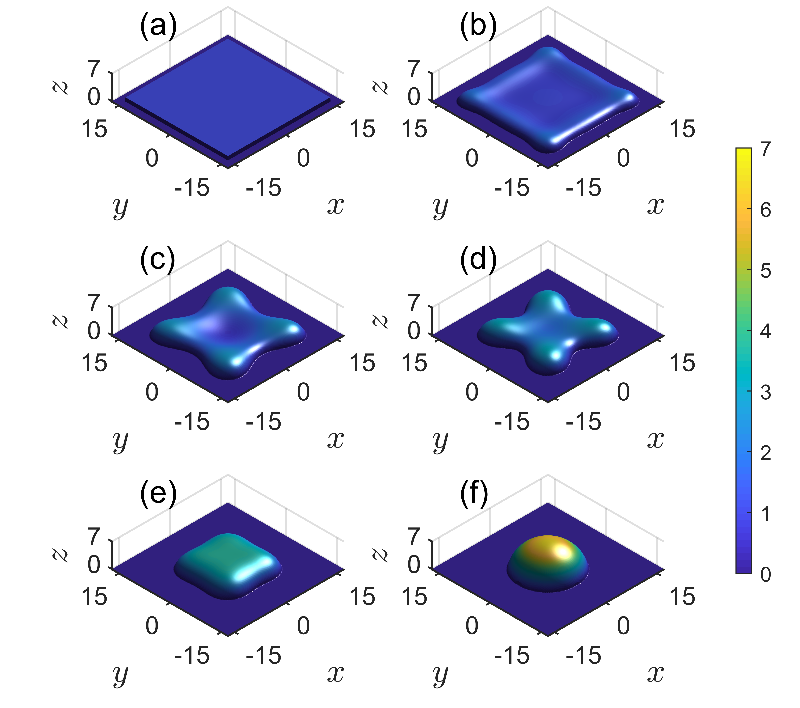}
\vspace*{-3mm}
\caption{Snapshots of the evolution of an initially \((30,30,1)\) cuboid island toward its equilibrium shape: (a) \(t=0\), (b) \(t=10\), (c) \(t=50\), (d) \(t=100\), (e) \(t=200\), (f) \(t=500\). The uniform mesh consists of 231226 triangles and 116295 vertices, and the time step is \(\tau=0.1\).}
\label{fig:square_3030}
\end{figure}

\begin{figure}[htp!]
\centering
\hspace*{-2mm}\includegraphics[width=.5\textwidth]{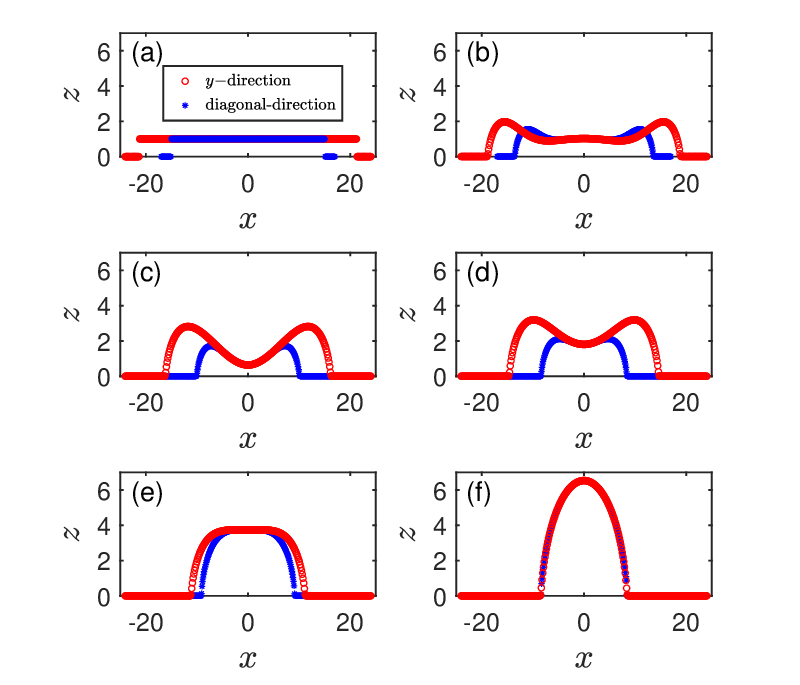}
\vspace*{-8mm}
\caption{Cross-sections of the island film along the \(y\)-direction and the diagonal direction for the example shown in Fig.~\ref{fig:square_3030}.}
\label{fig:square_3030_1D}
\end{figure}

As observed experimentally \cite{Thompson12,Ye10,Ye11b} and reported numerically in \cite{Jiang12,Naffouti17,Jiang20}, the retracting corners of a square island typically lag behind its edges during the early stage of dewetting. Fig.~\ref{fig:square_3030} shows the evolution of a relatively small square island of size \((30,30,1)\). At early times, mass rapidly accumulates near the corners, see Fig.~\ref{fig:square_3030}(b)--\ref{fig:square_3030}(d). As the evolution proceeds, the corners gradually catch up with the edges, the contact line becomes nearly circular, and the island eventually relaxes to a single cap-shaped equilibrium, see Fig.~\ref{fig:square_3030}(e)--\ref{fig:square_3030}(f). To further illustrate this evolution, the corresponding cross-sections in the \(y\)-direction and along the diagonal are shown in Fig.~\ref{fig:square_3030_1D}.

\begin{figure}[htp!]
\centering
\includegraphics[width=.5\textwidth]{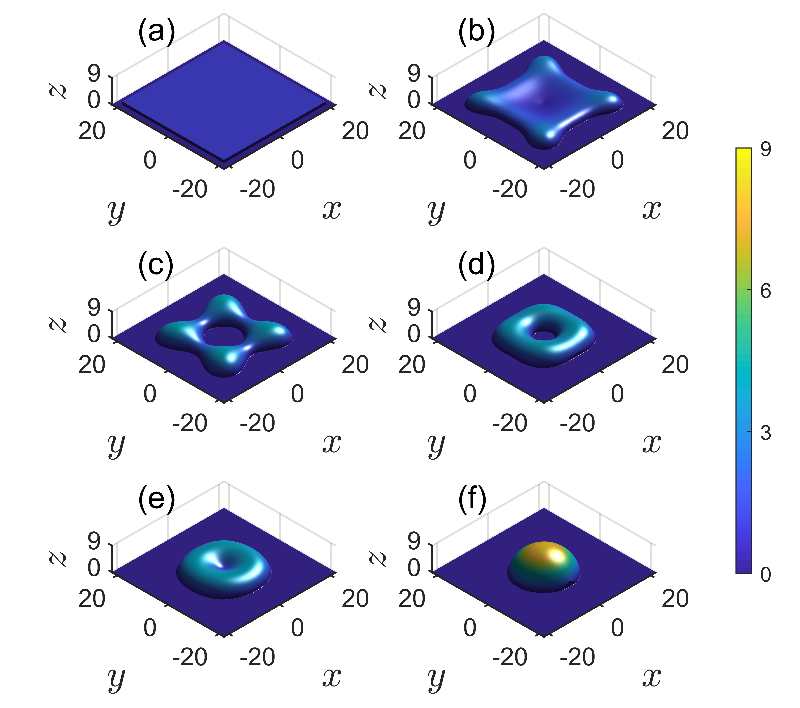}
\vspace*{-6mm}
\caption{Snapshots of the evolution of an initially \((40,40,1)\) cuboid island toward its equilibrium shape: (a) \(t=0\), (b) \(t=85\), (c) \(t=200\), (d) \(t=600\), (e) \(t=729\), (f) \(t=3000\). The uniform mesh consists of 387243 triangles and 194504 vertices, and the time step is \(\tau=0.1\).}
\label{fig:square_4040}
\end{figure}

\begin{figure}[htp!]
\centering
\includegraphics[width=.5\textwidth]{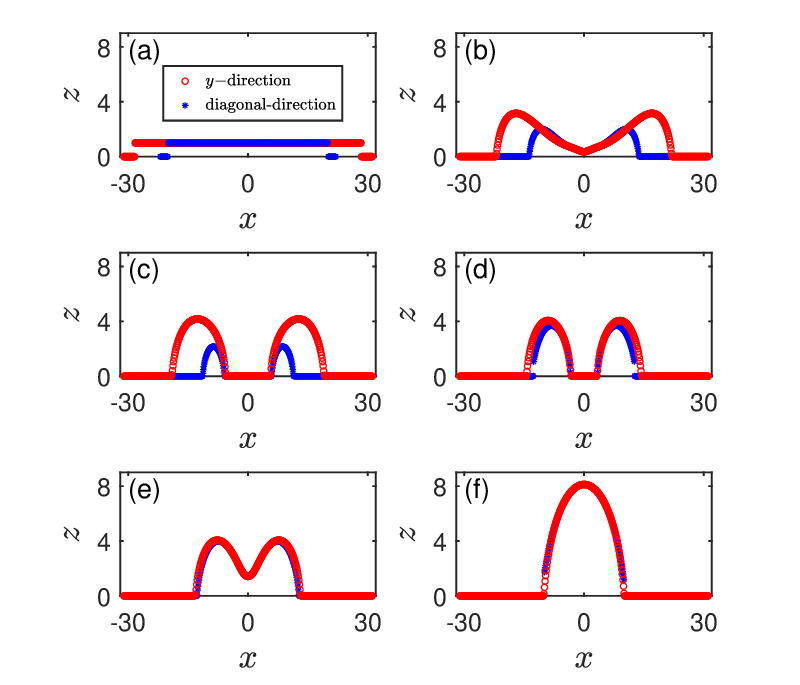}
\vspace*{-8mm}
\caption{Cross-sections of the island film along the \(y\)-direction and the diagonal direction for the example shown in Fig.~\ref{fig:square_4040}.}
\label{fig:square_4040_1D}
\end{figure}

We next increase the island size to \((40,40,1)\). In this case, a valley rapidly develops at the center and continues to deepen until it reaches the substrate, leading to hole formation, see Fig.~\ref{fig:square_4040}. As the evolution proceeds, the resulting ring-like structure becomes nearly axisymmetric, shrinks inward, and eventually merges into a single island. The inward shrinkage of such toroidal structures has also been reported in \cite{Zhao19,Jiang19a}. The corresponding cross-sections are shown in Fig.~\ref{fig:square_4040_1D}.

\begin{figure}[htp!]
\centering
\hspace*{-3mm}\includegraphics[width=.5\textwidth]{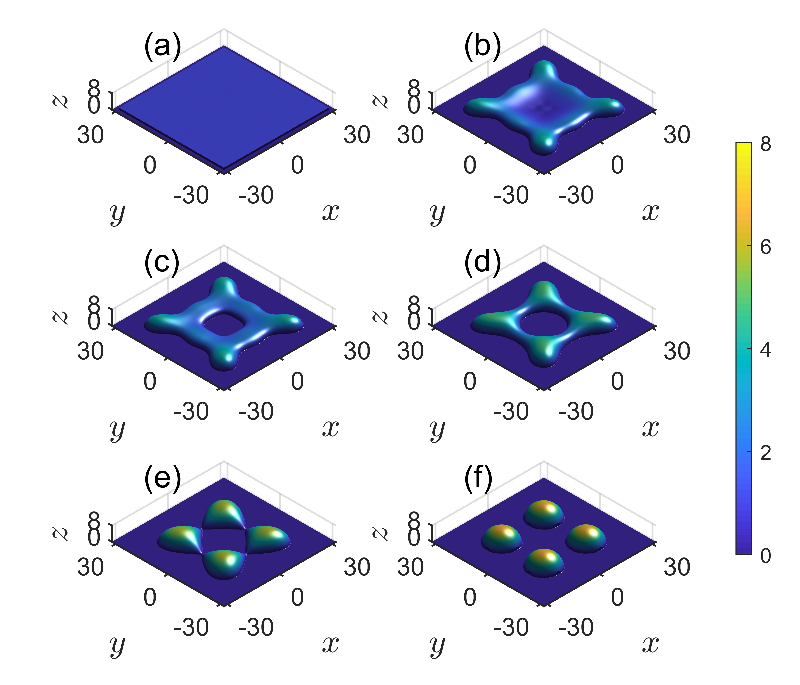}
\vspace*{-7mm}
\caption{Snapshots of the evolution of an initially \((60,60,1)\) cuboid island toward its equilibrium shape: (a) \(t=0\), (b) \(t=260\), (c) \(t=280\), (d) \(t=500\), (e) \(t=800\), (f) \(t=3000\).}
\label{fig:square_6060}
\end{figure}

\begin{figure}[htp!]
\centering
\hspace*{-1mm}\includegraphics[width=.5\textwidth]{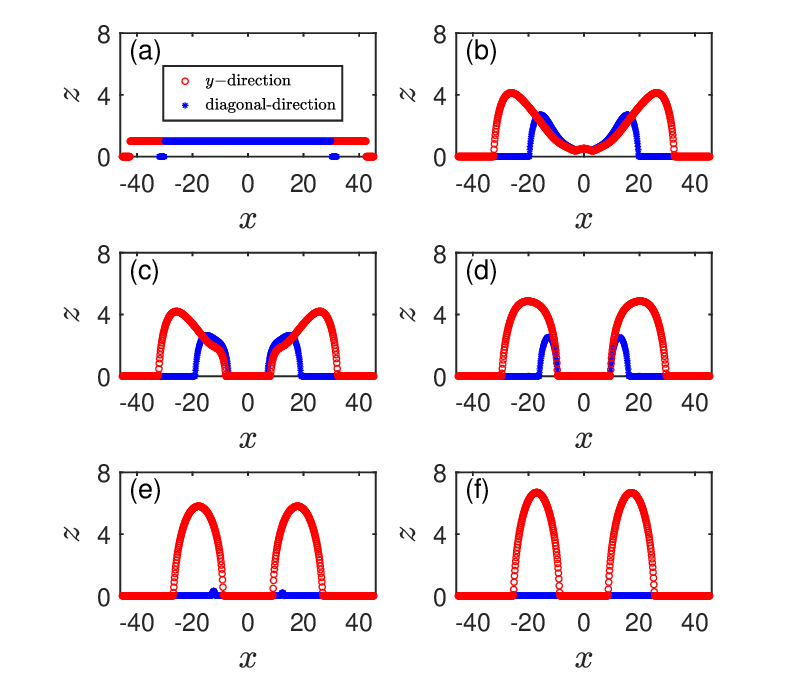}
\vspace*{-8mm}
\caption{Cross-sections of the island film along the \(y\)-direction and the diagonal direction for the example shown in Fig.~\ref{fig:square_6060}.}
\label{fig:square_6060_1D}
\end{figure}

For an even larger square island, namely \((60,60,1)\), the film again develops a central hole, but the subsequent evolution is markedly different. Instead of forming an approximately axisymmetric ring that shrinks into a single island, the film breaks into four isolated islands, see Fig.~\ref{fig:square_6060}(b)--\ref{fig:square_6060}(e). These islands then continue to relax toward their equilibrium shapes, see Fig.~\ref{fig:square_6060}(f). The corresponding cross-sections are shown in Fig.~\ref{fig:square_6060_1D}.

These results indicate a clear size effect in the three-dimensional evolution of square islands. As the island size increases, the morphology changes from relaxation toward a single equilibrium island to hole formation and eventual breakup into several isolated islands.

\subsection{Evolution of elongated cuboid islands}
\begin{figure}[htp!]
\centering
\includegraphics[width=.5\textwidth]{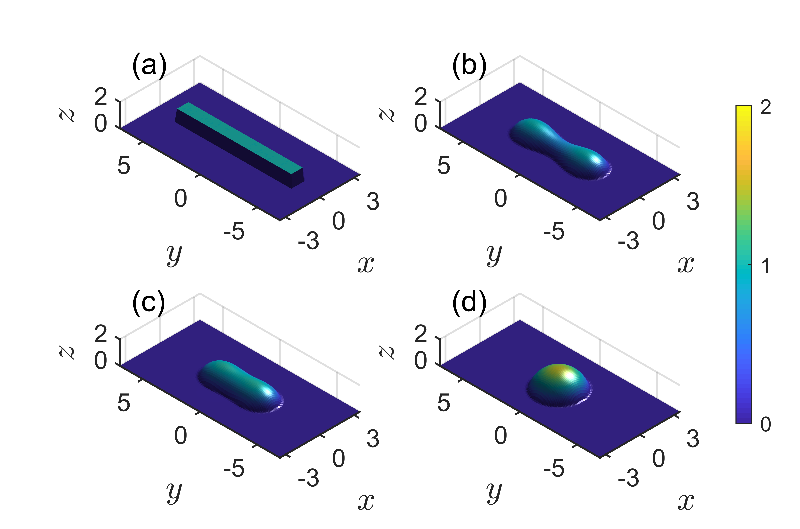}
\vspace*{-5mm}
\caption{Snapshots of the evolution of an initially \((1,10,1)\) cuboid island toward its equilibrium shape. (a) \(t=0\), (b) \(t=1\), (c) \(t=2\), (d) \(t=50\).}
\label{fig:long_10}
\end{figure}

We next consider elongated cuboid islands to examine the effect of increasing aspect ratio. For a relatively short cuboid island of size \((1,10,1)\), the film retracts rapidly and mass accumulates near the two ends, as shown in Fig.~\ref{fig:long_10}. As the evolution proceeds, the accumulated material moves toward the center, and the island eventually relaxes to a single cap-shaped equilibrium.

\begin{figure}[htp!]
\centering
\includegraphics[width=.5\textwidth]{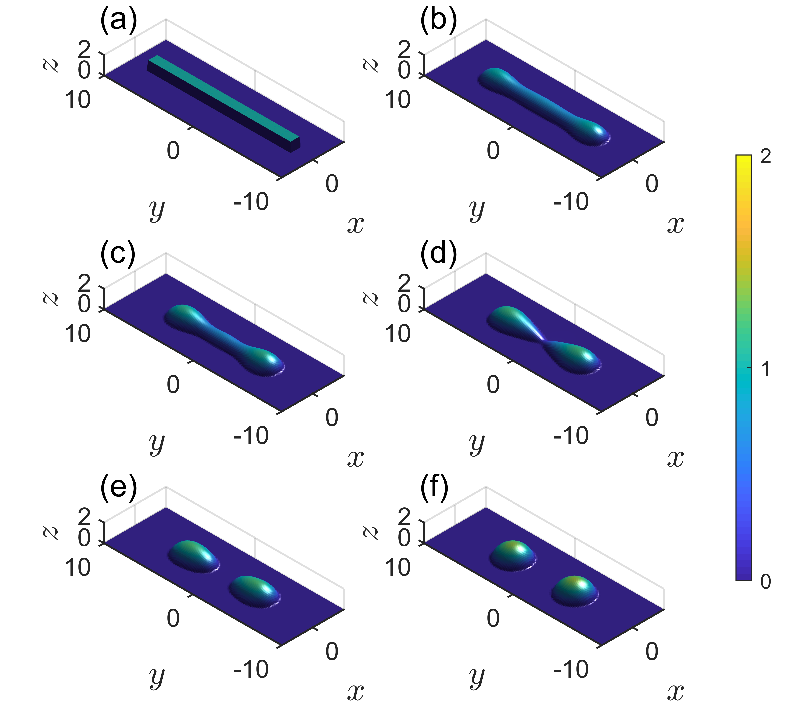}
\vspace*{-5mm}
\caption{Snapshots of the evolution of an initially \((1,16,1)\) cuboid island toward its equilibrium shape. (a) \(t=0\), (b) \(t=1\), (c) \(t=2\), (d) \(t=2.76\), (e) \(t=3\), (f) \(t=30\).}
\label{fig:long_16}
\end{figure}

\begin{figure}[htp!]
\centering
\includegraphics[width=.5\textwidth]{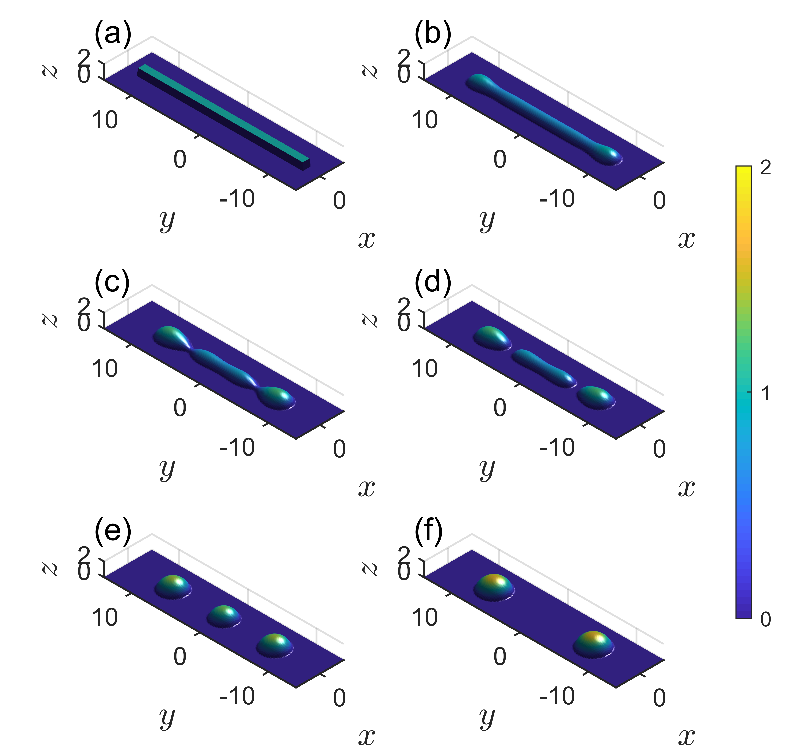}
\vspace*{-5mm}
\caption{Snapshots of the evolution of an initially \((1,24,1)\) cuboid island toward its equilibrium shape. (a) \(t=0\), (b) \(t=2\), (c) \(t=2.9\), (d) \(t=3\), (e) \(t=10\), (f) \(t=30\).}
\label{fig:long_24}
\end{figure}

As the cuboid becomes longer, pinch-off occurs and the film breaks into several isolated particles, as shown in Figs.~\ref{fig:long_16} and \ref{fig:long_24}. The breakup in the transverse direction can be viewed as a Rayleigh-like instability \cite{Kim15,Rayleigh1878}. For the case \((1,16,1)\), the film splits into two particles, which then relax toward their equilibrium shapes, as shown in Fig.~\ref{fig:long_16}. When the length is further increased to \((1,24,1)\), more particles are formed, and a clear coarsening process is observed, in which larger particles gradually absorb smaller ones, as shown in Fig.~\ref{fig:long_24}.

These examples show a clear aspect-ratio effect in the evolution of elongated cuboid islands. As the length increases, the morphology changes from retraction toward a single island to pinch-off, particle formation, and coarsening.

\subsection{Evolution of complex-shaped islands}

We finally consider more complex initial geometries, including square-ring and cross-shaped islands, in order to further demonstrate the capability of the proposed model in capturing rich three-dimensional dewetting morphologies.

\begin{figure}[htp!]
\centering
\hspace*{-2mm}\includegraphics[width = .5\textwidth]{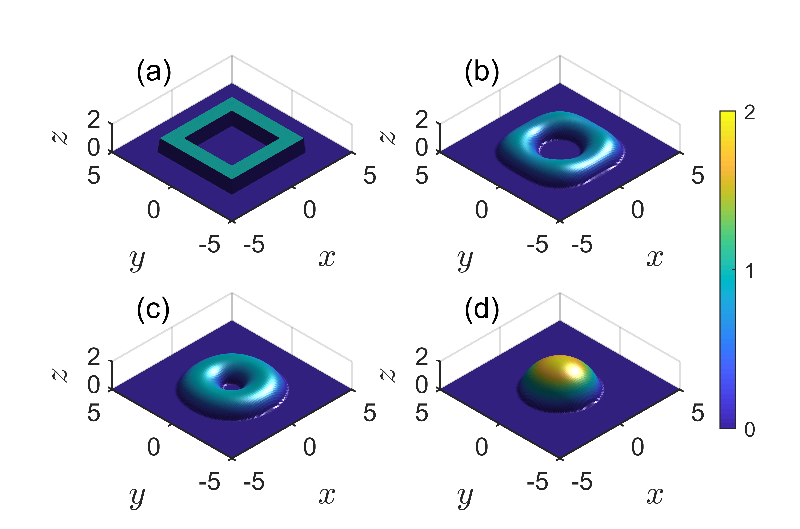}
\vspace*{-5mm}
\caption{Snapshots of the evolution of an initially square-ring island obtained from a \((6,6,1)\) cuboid by removing a \((4,4,1)\) cuboid from the center: (a) \(t=0\), (b) \(t=1\), (c) \(t=2\), (d) \(t=30\).}
\label{fig:ring_5}
\end{figure}

\begin{figure}[htp!]
\centering
\hspace*{-1mm}\includegraphics[width = .5\textwidth]{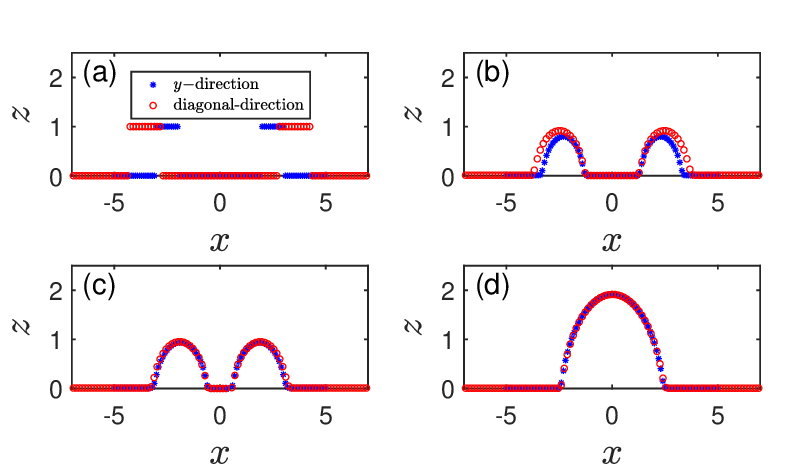}
\vspace*{-5mm}
\caption{Cross-sections of the island geometry for the example shown in Fig.~\ref{fig:ring_5}.}
\label{fig:ring_5_1D}
\end{figure}

We first study square-ring islands. The initial geometry is constructed from a square island of size \((c,c,1)\) by removing a concentric square island of size \((c-2,c-2,1)\), where \(c>0\). For a relatively small square-ring island with \(c=6\), the film evolves rapidly, and the ring gradually becomes nearly axisymmetric, shrinks inward, and finally merges into a single island, as shown in Fig.~\ref{fig:ring_5}. This behavior is similar to that observed for the square island in Fig.~\ref{fig:square_4040}. The corresponding cross-sections are shown in Fig.~\ref{fig:ring_5_1D}.

\begin{figure}[htp!]
\centering
\hspace*{-2mm}\includegraphics[width = .49\textwidth]{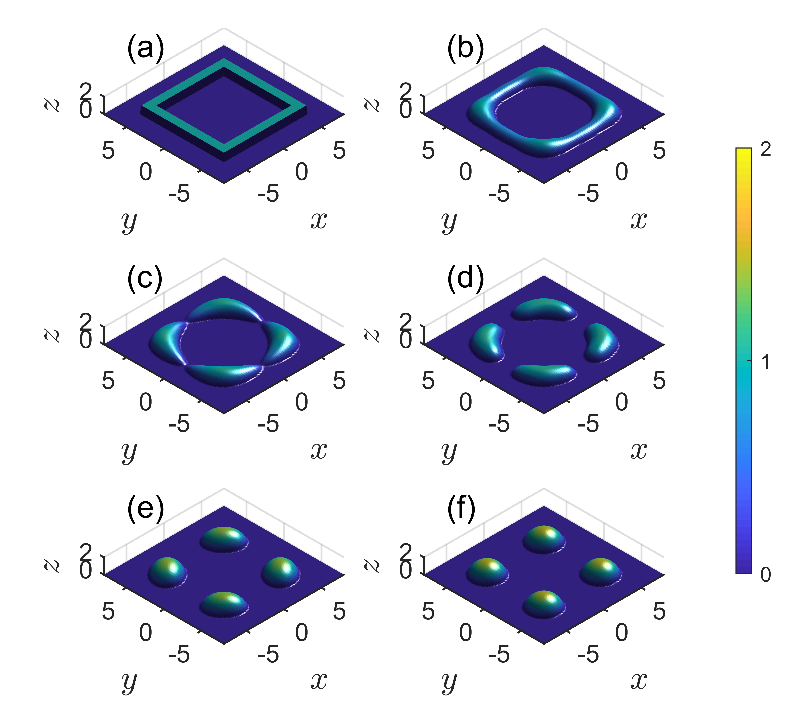}
\vspace*{-5mm}
\caption{Snapshots of the evolution of an initially square-ring island obtained from a \((11,11,1)\) cuboid by removing a \((9,9,1)\) cuboid from the center: (a) \(t=0\), (b) \(t=1\), (c) \(t=1.79\), (d) \(t=2\), (e) \(t=3\), (f) \(t=10\).}
\label{fig:ring_10}
\end{figure}

\begin{figure}[htp!]
\centering
\hspace*{-2mm}\includegraphics[width = .47\textwidth]{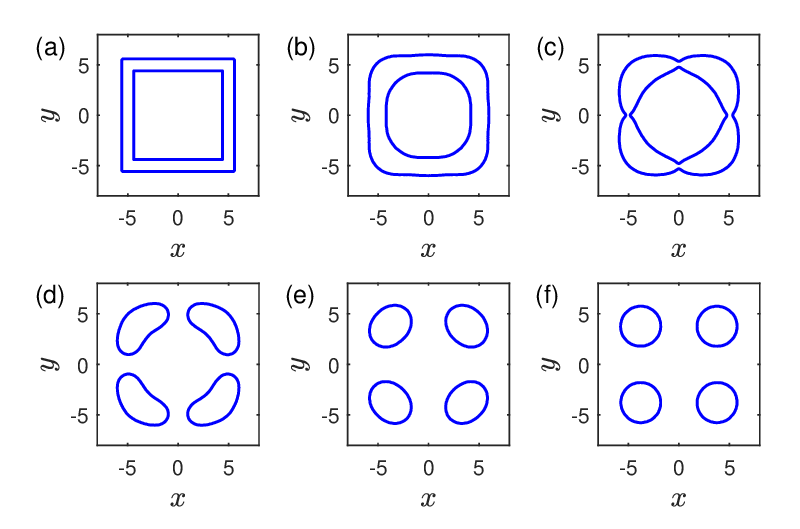}
\vspace*{-3mm}
\caption{Cross-sections of the island geometry for the example shown in Fig.~\ref{fig:ring_10}.}
\label{fig:ring_10_1D}
\end{figure}

\begin{figure}[htp!]
\centering
\hspace*{-2mm}\includegraphics[width = .5\textwidth]{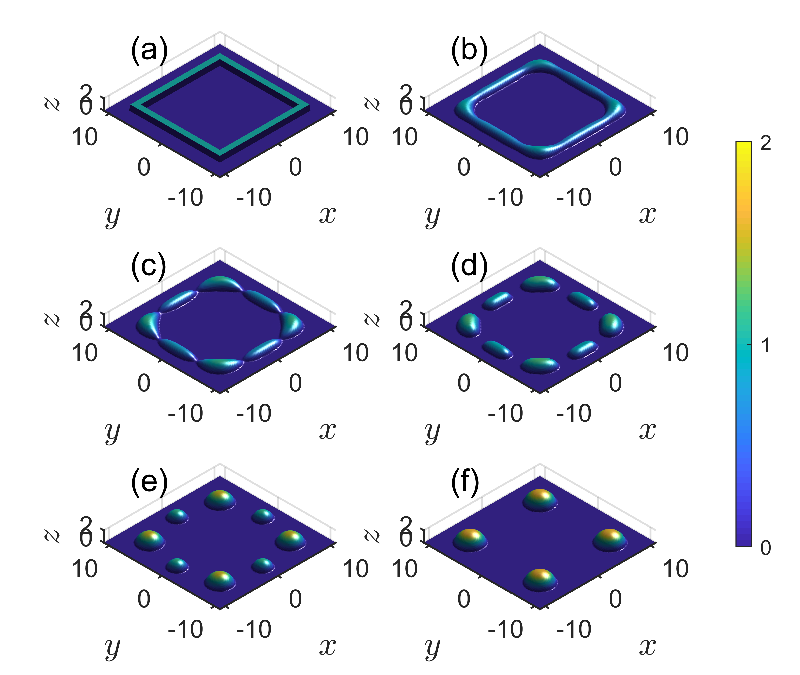}
\vspace*{-5mm}
\caption{Snapshots of the evolution of an initially square-ring island obtained from a \((17,17,1)\) cuboid by removing a \((15,15,1)\) cuboid from the center: (a) \(t=0\), (b) \(t=1\), (c) \(t=2.15\), (d) \(t=2.5\), (e) \(t=5\), (f) \(t=30\).}
\label{fig:ring_16}
\end{figure}

\begin{figure}[htp!]
\centering
\hspace*{-2mm}\includegraphics[width = .5\textwidth]{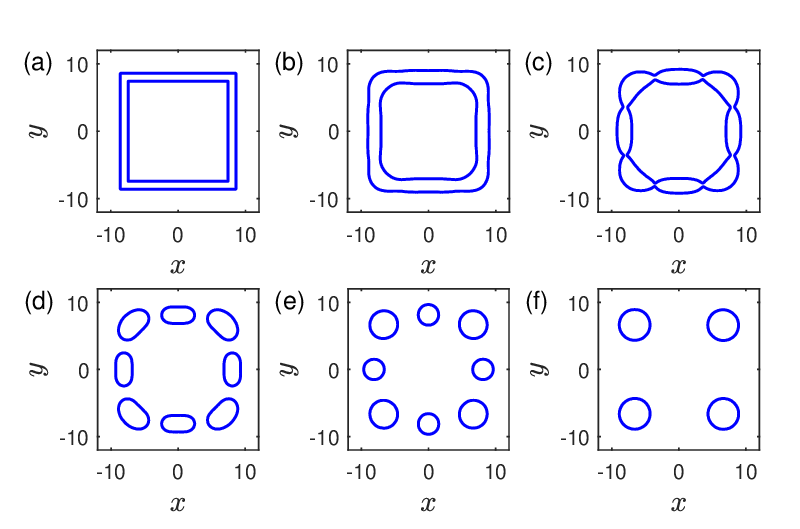}
\vspace*{-5mm}
\caption{Cross-sections of the island geometry for the example shown in Fig.~\ref{fig:ring_16}.}
\label{fig:ring_16_1D}
\end{figure}

As the size of the square-ring island increases, more complicated topological changes occur. For \(c=11\), pinch-off takes place and the film breaks into four isolated islands, which then evolve toward their equilibrium shapes, as shown in Fig.~\ref{fig:ring_10}. When the size is further increased to \(c=17\), more isolated islands are formed and a clear coarsening process is observed, as shown in Fig.~\ref{fig:ring_16}. To better illustrate these evolutions, the corresponding cross-sections are plotted in Fig.~\ref{fig:ring_10_1D} and Fig.~\ref{fig:ring_16_1D}.

\begin{figure}[htp!]
\centering
\includegraphics[width = .5\textwidth]{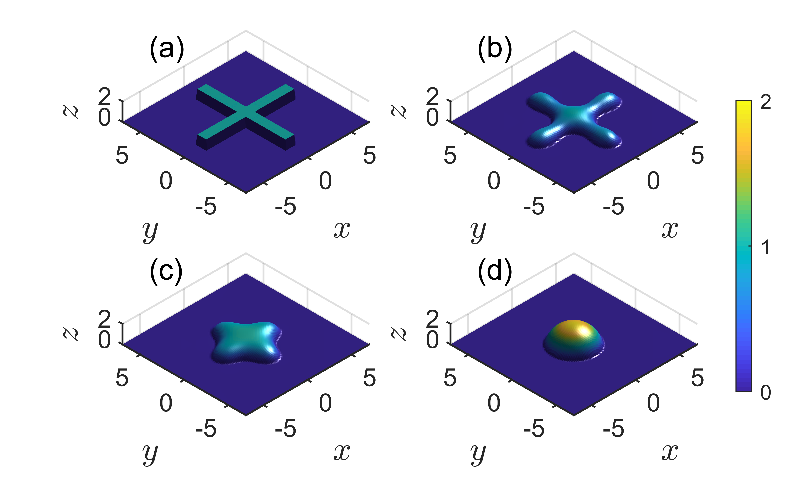}
\vspace*{-5mm}
\caption{Snapshots of the evolution of an initially cross-shaped island toward its equilibrium shape, where the initial island consists of four \((1,4,1)\) cuboids forming the limbs and one \((1,1,1)\) cube at the center: (a) \(t=0\), (b) \(t=0.2\), (c) \(t=1\), (d) \(t=20\).}
\label{fig:shi_5}
\end{figure}

\begin{figure}[htp!]
\centering
\includegraphics[width = .5\textwidth]{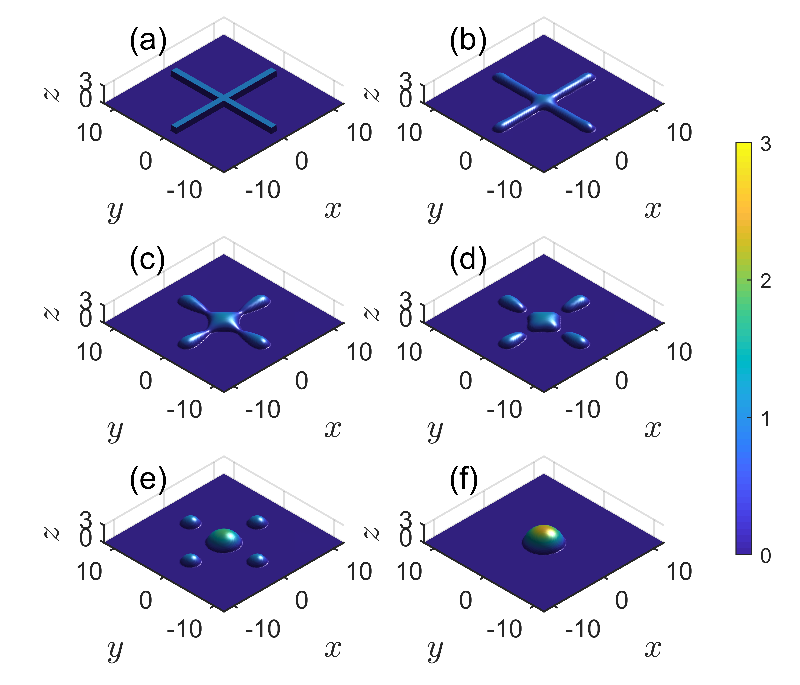}
\vspace*{-5mm}
\caption{Snapshots of the evolution of an initially cross-shaped island toward its equilibrium shape, where the initial island consists of four \((1,9,1)\) cuboids forming the limbs and one \((1,1,1)\) cube at the center: (a) \(t=0\), (b) \(t=0.2\), (c) \(t=1.05\), (d) \(t=1.2\), (e) \(t=5\), (f) \(t=50\).}
\label{fig:shi_10}
\end{figure}

We next consider cross-shaped islands. The initial geometry consists of a unit cube together with four identical \((1,c,1)\) cuboids as its limbs. For \(c=4\), the four limbs retract quickly and the cross-shaped island eventually evolves into a single island, as shown in Fig.~\ref{fig:shi_5}. When the limb length is increased to \(c=9\), the island breaks into five isolated particles. As time evolves, the largest particle in the center gradually absorbs the smaller ones, indicating a coarsening process, as shown in Fig.~\ref{fig:shi_10}.

These examples indicate that the proposed model can also handle complex initial geometries and capture rich three-dimensional topological changes, including inward contraction, pinch-off, the formation of multiple isolated islands, and coarsening.

\section{Conclusions}

In this paper, we proposed a sharp-interface model with wetting potential for describing solid-state dewetting dynamics from an energy variational perspective. To solve the resulting model efficiently, we developed a finite element method with a semi-implicit time discretization. Numerical experiments in both two and three dimensions demonstrated that the proposed method is effective and capable of capturing essential features of the dewetting process.

With the proposed model, a variety of complex phenomena were successfully simulated, including pinch-off, Rayleigh-like instability, edge retraction, and corner mass accumulation. Relative to the phase-field model \cite{Jiang12}, the present approach is computationally more efficient since it reduces the spatial dimension by one. In addition, unlike thickness-independent sharp-interface models \cite{Wang15,Bao17a,Jiang20,Zhao20}, the proposed model can naturally capture topological changes and handle complex geometries.

Several issues deserve further investigation. One is to extend the present model and numerical method to the anisotropic case. Another is to obtain a better understanding of the asymptotic behavior as \(\varepsilon \to 0\). In addition, for large-scale dewetting problems, the computational cost remains high, especially in three dimensions, and more efficient numerical strategies are still needed.

\begin{acknowledgments}
W.H. acknowledges support from the National Natural Science Foundation of China
under Grant No. 12001034. X.R. acknowledges support from the National Natural
Science Foundation of China under Grant No. 12201436.
\end{acknowledgments}

\appendix

\section{First variation of the free energy in two dimensions}\label{Appendix:A}

In this appendix, we derive the first variation of the energy functional in Eq.~\eqref{eq:energy_functional} with respect to the height function \(h(x)\). For any \(g\in H^1(I)\) with \(I=(a,b)\), integration by parts yields
\begin{equation}
\begin{aligned}
& \left.\frac{d}{d\alpha} W^\varepsilon(h+\alpha g)\right|_{\alpha=0}\\
=& \int_a^b \left[(\gamma^\varepsilon)'(h)\sqrt{1+(\partial_x h)^2}\,g
+ \frac{\gamma^\varepsilon(h)\,\partial_x h\,\partial_x g}{\sqrt{1+(\partial_x h)^2}}\right]\,dx \\
=& \int_a^b \left[(\gamma^\varepsilon)'(h)\sqrt{1+(\partial_x h)^2}
-\partial_x\!\left(\frac{\gamma^\varepsilon(h)\,\partial_x h}{\sqrt{1+(\partial_x h)^2}}\right)\right] g\,dx \\
&\quad + \frac{\gamma^\varepsilon(h)\,\partial_x h}{\sqrt{1+(\partial_x h)^2}}\,g\Big|_{x=a}^{x=b} \\
=& \int_a^b \left[
\frac{(\gamma^\varepsilon)'(h)}{\sqrt{1+(\partial_x h)^2}}
-\frac{\gamma^\varepsilon(h)\,\partial_{xx}h}{\bigl(1+(\partial_x h)^2\bigr)^{3/2}}
\right] g\,dx \\
&\quad + \frac{\gamma^\varepsilon(h)\,\partial_x h}{\sqrt{1+(\partial_x h)^2}}\,g\Big|_{x=a}^{x=b}.
\end{aligned}
\end{equation}
Using the curvature formula
 \[
\kappa=-\frac{\partial_{xx}h}{\bigl(1+(\partial_x h)^2\bigr)^{3/2}},
\]
we obtain the variational derivative in \eqref{eq:chemical_potential}. Under the boundary conditions \eqref{bc:h_x}, the boundary term vanishes.

\section{Mass conservation and energy dissipation}\label{Appendix:B}

\paragraph*{Mass conservation.}
The total mass of the thin film is defined by
\begin{equation}\label{eq:mass}
A(t)=\int_a^b h(x,t)\,dx.
\end{equation}
Then, by \eqref{eq:mass} and the first equation in \eqref{eq:governing_equation}, we have
\begin{equation}
\begin{aligned}
\frac{dA(t)}{dt}
=& \int_a^b \partial_t h(x,t)\,dx \\
=& \int_a^b \partial_x\left(\frac{\partial_x\mu}{\sqrt{1+(\partial_x h)^2}}\right)\,dx \\
=& \frac{\partial_x\mu}{\sqrt{1+(\partial_x h)^2}}\Big|_{x=a}^{x=b} \\
=& 0,
\end{aligned}
\end{equation}
under the boundary conditions \eqref{bc:mass_zeros}.

\paragraph*{Energy dissipation.}
The total free energy of the system is given by
\begin{equation}
W^\varepsilon(t)=\int_a^b \gamma^\varepsilon(h)\sqrt{1+(\partial_x h)^2}\,dx.
\end{equation}
Differentiating with respect to time gives
\begin{equation}
\begin{aligned}
\frac{dW^\varepsilon}{dt}
=& \int_a^b \left[(\gamma^\varepsilon)'(h)\,\partial_t h\,\sqrt{1+(\partial_x h)^2}
+ \frac{\gamma^\varepsilon(h)\,\partial_x h\,\partial_{xt}h}{\sqrt{1+(\partial_x h)^2}}\right]\,dx.
\end{aligned}
\end{equation}
Integrating by parts and using the boundary conditions \eqref{bc:h_x}, we obtain
\begin{equation}
\begin{aligned}
\frac{dW^\varepsilon}{dt}
=& \int_a^b \left[
\frac{(\gamma^\varepsilon)'(h)}{\sqrt{1+(\partial_x h)^2}}
-\frac{\gamma^\varepsilon(h)\,\partial_{xx}h}{\bigl(1+(\partial_x h)^2\bigr)^{3/2}}
\right]\partial_t h\,dx \\
=& \int_a^b \mu\,\partial_t h\,dx.
\end{aligned}
\end{equation}
Substituting the evolution equation in \eqref{eq:governing_equation} yields
\begin{equation}
\begin{aligned}
\frac{dW^\varepsilon}{dt}
=& \int_a^b \mu\,\partial_x\left(\frac{\partial_x\mu}{\sqrt{1+(\partial_x h)^2}}\right)\,dx \\
=& -\int_a^b \frac{1}{\sqrt{1+(\partial_x h)^2}}\,(\partial_x\mu)^2\,dx \\
\le& 0,
\end{aligned}
\end{equation}
where the boundary term vanishes due to \eqref{bc:mass_zeros}.

\section{Construction of \texorpdfstring{$\zeta(h)$}{zeta(h)}}\label{Appendix:C}

In this appendix, we construct the quadratic approximation \(\zeta(h)\) used in the semi-implicit scheme.
Recall that
\begin{equation}
(\gamma^\varepsilon)'(h)
= \frac{1-\sigma}{\varepsilon}\left[-e^{-h/\varepsilon}+e^{-h/(2\varepsilon)}\right].
\end{equation}
For \(h\le \bar h\ll 1\), we approximate \((\gamma^\varepsilon)'(h)\) by a quadratic polynomial of the form
\[
\zeta(h)=c_1h+c_2h^2,
\]
which automatically satisfies \(\zeta(0)=0=(\gamma^\varepsilon)'(0)\). The coefficients \(c_1\) and \(c_2\) are chosen so that
\begin{equation}\label{eq:zeta_condition}
\zeta(\bar h)=(\gamma^\varepsilon)'(\bar h),\qquad
\zeta'(\bar h)=(\gamma^\varepsilon)''(\bar h).
\end{equation}

Since \((\gamma^\varepsilon)'(0)=0\), it is convenient to introduce
\[
f(h)=\frac{(\gamma^\varepsilon)'(h)}{h}
=\frac{1-\sigma}{\varepsilon}\,
\frac{-e^{-h/\varepsilon}+e^{-h/(2\varepsilon)}}{h}.
\]
Then
\[
(\gamma^\varepsilon)'(h)=h\,f(h).
\]
A first-order Taylor expansion of \(f(h)\) at \(h=\bar h\) gives
\begin{equation}
f(h)=f(\bar h)+f'(\bar h)(h-\bar h)+o(h-\bar h).
\end{equation}
Motivated by this expansion, we define
\begin{equation}
\begin{aligned}
\zeta(h)
=&\, h\Bigl[f(\bar h)+f'(\bar h)(h-\bar h)\Bigr] \\
=&\, \bigl(f(\bar h)-\bar h f'(\bar h)\bigr)h + f'(\bar h)h^2 \\
\triangleq&\, c_1 h + c_2 h^2.
\end{aligned}
\end{equation}
Therefore,
\[
c_1=f(\bar h)-\bar h f'(\bar h),\qquad c_2=f'(\bar h).
\]
By direct computation, we obtain
\begin{equation}
\begin{aligned}
c_1
=&\, \frac{1-\sigma}{\varepsilon}\left[
2\frac{-e^{-\bar h/\varepsilon}+e^{-\bar h/(2\varepsilon)}}{\bar h}
-\frac{e^{-\bar h/\varepsilon}-\frac12 e^{-\bar h/(2\varepsilon)}}{\varepsilon}
\right], \\
c_2
=&\, \frac{1-\sigma}{\varepsilon}\left[
\frac{e^{-\bar h/\varepsilon}-\frac12 e^{-\bar h/(2\varepsilon)}}{\varepsilon \bar h}
+\frac{e^{-\bar h/\varepsilon}-e^{-\bar h/(2\varepsilon)}}{\bar h^2}
\right].
\end{aligned}
\end{equation}

We now verify that \(\zeta(h)\) satisfies \eqref{eq:zeta_condition}.
First, since \(\zeta(h)=c_1h+c_2h^2\), it follows immediately that
\[
\zeta(0)=0=(\gamma^\varepsilon)'(0).
\]
Next, evaluating \(\zeta(h)\) at \(h=\bar h\), we obtain
\[
\zeta(\bar h)=\bar h\,f(\bar h)=(\gamma^\varepsilon)'(\bar h).
\]
Moreover,
\[
\zeta'(h)=f(\bar h)-\bar h f'(\bar h)+2f'(\bar h)h,
\]
and hence
\[
\zeta'(\bar h)=f(\bar h)+\bar h f'(\bar h).
\]
Since \((\gamma^\varepsilon)'(h)=h\,f(h)\), we have
\[
(\gamma^\varepsilon)''(h)=f(h)+h f'(h),
\]
which implies
\[
\zeta'(\bar h)=(\gamma^\varepsilon)''(\bar h).
\]
Therefore, \(\zeta(h)\) satisfies all the conditions in \eqref{eq:zeta_condition}.

\section{First variation of the free energy in three dimensions}\label{Appendix:D}

Let \(\Omega=[a,b]\times[c,d]\) be a fixed domain, and let the film-vapor interface
be represented by the graph
\[
\Gamma=\{(x,y,h(x,y,t)):(x,y)\in\Omega\}.
\]
The total free energy of the system is
\begin{equation}
W^\varepsilon(h)=\int_\Omega \gamma^\varepsilon(h)\sqrt{1+|\nabla h|^2}\,dx\,dy.
\end{equation}
For any \(g\in H^1(\Omega)\), integrating by parts gives
\begin{equation}
\begin{aligned}
& \left.\frac{d}{d\alpha}W^\varepsilon(h+\alpha g)\right|_{\alpha=0}\\
=& \int_\Omega \left[
(\gamma^\varepsilon)'(h)\,g\,\sqrt{1+|\nabla h|^2}
+\gamma^\varepsilon(h)\frac{\nabla h\cdot \nabla g}{\sqrt{1+|\nabla h|^2}}
\right]\,dx\,dy \\
=& \int_\Omega (\gamma^\varepsilon)'(h)\sqrt{1+|\nabla h|^2}\,g
-\nabla\cdot\left(
\frac{\gamma^\varepsilon(h)\nabla h}{\sqrt{1+|\nabla h|^2}}
\right) g\,dx\,dy \\
&\quad + \int_{\partial\Omega}
\frac{\gamma^\varepsilon(h)\,\nabla h\cdot \mathbf n}{\sqrt{1+|\nabla h|^2}}\,g\,ds.
\end{aligned}
\end{equation}
Expanding the divergence term and using the mean curvature expression
\[
\kappa
=
-\frac{(1+(\partial_y h)^2)\partial_{xx}h
-2\partial_x h\,\partial_y h\,\partial_{xy}h
+(1+(\partial_x h)^2)\partial_{yy}h}
{(1+|\nabla h|^2)^{3/2}},
\]
we obtain
\begin{equation}
\begin{aligned}
\left.\frac{d}{d\alpha}W^\varepsilon(h+\alpha g)\right|_{\alpha=0}
=&
\int_\Omega \left[
\frac{(\gamma^\varepsilon)'(h)}{\sqrt{1+|\nabla h|^2}}
+\gamma^\varepsilon(h)\kappa
\right] g dxdy \\
&\quad +
\int_{\partial\Omega}
\frac{\gamma^\varepsilon(h)\,\nabla h\cdot \mathbf n}{\sqrt{1+|\nabla h|^2}}\,g\,ds.
\end{aligned}
\end{equation}
The boundary term vanishes after imposing the boundary condition \(\nabla h\cdot \mathbf n=0\) on \(\partial\Omega\). Hence the variational derivative of \(W^\varepsilon\) with respect to \(h\) is
\begin{equation}
\frac{\delta W^\varepsilon}{\delta h}
=
\gamma^\varepsilon(h)\kappa
+
\frac{(\gamma^\varepsilon)'(h)}{\sqrt{1+|\nabla h|^2}}.
\end{equation}
We therefore define the chemical potential by
\begin{equation}
\mu=
\gamma^\varepsilon(h)\kappa
+
\frac{(\gamma^\varepsilon)'(h)}{\sqrt{1+|\nabla h|^2}}.
\end{equation}
For solid-state dewetting, the interface evolution is governed by surface diffusion. With the surface flux \(\mathbf J_s=-\nabla_s\mu\), mass conservation on the surface gives \cite{BaoZhao23}
\[
v_n=-\nabla_s\cdot \mathbf J_s=\Delta_s\mu.
\]
On the other hand, for the graph representation \(X(x,y,t)=(x,y,h(x,y,t))\), the unit normal vector is
\[
\mathbf n=\frac{1}{\sqrt{1+|\nabla h|^2}}(-\partial_x h,-\partial_y h,1)^T,
\]
and thus
\begin{equation}
v_n=X_t\cdot \mathbf n=\frac{h_t}{\sqrt{1+|\nabla h|^2}}.
\end{equation}
Combining the above two relations, we arrive at
\begin{equation}
h_t=\sqrt{1+|\nabla h|^2}\,\Delta_s\mu.
\end{equation}
Therefore, the three-dimensional graph model reads
\begin{equation}
\left\{
\begin{aligned}
& h_t = \sqrt{1 + |\nabla h|^2}\,\Delta_s \mu,\\
& \mu = \gamma^\varepsilon(h)\kappa + \dfrac{(\gamma^\varepsilon)'(h)}{\sqrt{1 + |\nabla h|^2}},
\end{aligned}
\right.
\quad (x,y)\in\Omega,\quad t>0.
\end{equation}
Together with the boundary conditions
\begin{gather}
\nabla h\cdot \mathbf n = 0 \qquad \text{on } \partial\Omega,\\
\nabla_s\mu\cdot \mathbf c = 0 \qquad \text{on } \partial S,
\end{gather}
where \(S=\{(x,y,h(x,y,t)):(x,y)\in\Omega\}\) denotes the film-vapor interface, and \(\mathbf c\) is the outward co-normal vector along \(\partial S\).

\bibliography{mybib}

\end{document}